\documentclass[letterpaper, 11pt]{amsart}

\usepackage{amsmath,amsthm,amsfonts,amssymb,amscd}
\usepackage{bbm}
\usepackage{bm}
\usepackage{tikz}
\usepackage{tikz-cd}
\usepackage{BOONDOX-calo}
\usepackage{standalone}
\usepackage[hidelinks]{hyperref}
\usetikzlibrary{arrows,chains,matrix,positioning,scopes}
\usepackage[letterpaper, left=3cm,right=3cm, top=3cm, bottom=3cm]{geometry}

\newcommand{\CC}{{\mathbb C}}

\newcommand{\ZZ}{{\mathbb Z}}
\newcommand{\QQ}{{\mathbb Q}}

\newcommand{\inEnd}{{\underline{End}}}
\newcommand{\fu}{{\mathfrak{u}}}
\newcommand{\fg}{{\mathfrak{g}}}
\newcommand{\inHom}{{\underline{Hom}}}

\newcommand\dual{\raise0.9ex\hbox{$\scriptscriptstyle\vee$}}

\theoremstyle{definition}

\theoremstyle{remark}
\newtheorem*{rem}{Remark}

\theoremstyle{plain}
\newtheorem{thm}{Theorem} 
\newtheorem{prop}[thm]{Proposition}
\newtheorem{lemma}[thm]{Lemma}  
\newtheorem{cor}[thm]{Corollary}
\numberwithin{thm}{subsection}

\newtheorem{manualtheoreminner}{Theorem}
\newenvironment{thm'}[1]{%
  \renewcommand\themanualtheoreminner{#1}%
  \manualtheoreminner
}{\endmanualtheoreminner}

\tikzset{>=stealth}

\makeatletter
\def\@seccntformat#1{%
  \protect\textup{\protect\@secnumfont
    \ifnum\pdfstrcmp{subsection}{#1}=0 \bfseries\fi
    \csname the#1\endcsname
    \protect\@secnumpunct
  }%
}  
\makeatother

\makeatletter
\@namedef{subjclassname@2020}{%
  $2020$ Mathematics Subject Classification}
\makeatother

\begin{document}

\title[On unipotent radicals of motivic Galois groups]{On unipotent radicals of motivic Galois groups}
\author{Payman Eskandari and V. Kumar Murty}
\address{Department of Mathematics, University of Toronto, 40 St. George St., Room 6290, Toronto, Ontario, Canada, M5S 2E4}
\email{payman@math.toronto.edu , murty@math.toronto.edu}
\begin{abstract}
Let $\mathbf{T}$ be a neutral Tannakian category over a field of characteristic zero with unit object $\mathbbm{1}$, and equipped with a filtration $W_\cdot$ similar to the weight filtration on mixed motives. Let $M$ be an object of $\mathbf{T}$, and $\underline{\fu}(M)\subset W_{-1}\inHom(M,M)$ the Lie algebra of the kernel of the natural surjection from the fundamental group of $M$ to the fundamental group of $Gr^WM$. A result of Deligne gives a characterization of $\underline{\fu}(M)$ in terms of the extensions $0\longrightarrow W_pM \longrightarrow M \longrightarrow M/W_pM \longrightarrow 0$: it states that $\underline{\fu}(M)$ is the smallest subobject of $W_{-1}\inHom(M,M)$ such that the sum of the aforementioned extensions, considered as extensions of $\mathbbm{1}$ by $W_{-1}\inHom(M,M)$, is the pushforward of an extension of $\mathbbm{1}$ by $\underline{\fu}(M)$. In this article, we study each of the above-mentioned extensions individually in relation to $\underline{\fu}(M)$. Among other things, we obtain a refinement of Deligne's result, where we give a sufficient condition for when an individual extension $0\longrightarrow W_pM \longrightarrow M \longrightarrow M/W_pM \longrightarrow 0$ is the pushforward of an extension of $\mathbbm{1}$ by $\underline{\fu}(M)$. In the second half of the paper, we give an application to mixed motives whose unipotent radical of the motivic Galois group is as large as possible (i.e. with $\underline{\fu}(M)= W_{-1}\inHom(M,M)$). Using Grothedieck's formalism of {\it extensions panach\'{e}es} we prove a classification result for such motives. Specializing to the category of mixed Tate motives we obtain a classification result for 3-dimensional mixed Tate motives over $\QQ$ with three weights and large unipotent radicals.  
\end{abstract}
\maketitle

\section{Introduction}
\subsection{About this paper}
Let $\mathbf{T}$ be a neutral Tannakian category over a field $K$ of characteristic zero, equipped with a weight filtration $W_\cdot$ similar to the weight filtration on mixed motives (functorial, increasing, finite on every object, exact, and respecting the tensor structure). For 
example, one might keep in mind the category of mixed Hodge structures. In fact, this is a concrete example that illustrates well the 
main results.
\medskip\par
Let $M$ be an object of $\mathbf{T}$, and $\underline{\fu}(M)$ the Lie algebra of the kernel of the natural map from the fundamental group of $M$ to that of $Gr^WM$. A result of Deligne describes $\underline{\fu}(M)$ in terms of extensions that arise naturally from the weight filtration of $M$.
For each integer $p$, let $\mathcal{E}_p(M)$ be the extension
\begin{equation}\label{eq1into}
0 \ \longrightarrow \ W_pM   \ \longrightarrow \ M  \ \longrightarrow \ M/W_pM  \ \longrightarrow \ 0,
\end{equation}
considered as an element in $Ext^1(\mathbbm{1}, W_{-1}\inEnd(M))$ (where $\inEnd(M)$ means $\inHom(M,M)$, 
the latter being the internal Hom). Deligne characterizes $\underline{\fu}(M)$ in terms of the sum
\[
\mathcal{E}(M) \ := \ \sum\limits_p \mathcal{E}_p(M) \ \in \ Ext^1(\mathbbm{1}, W_{-1}\inEnd(M)).
\]
The first half of this paper refines this by  developing conditions under which the individual extensions $\mathcal{E}_p(M)$ can be
related to $\underline{\fu}(M)$.
\medskip\par
The second half of the paper specializes to the setting of mixed motives and gives an application of the first half to mixed motives whose unipotent radical of the motivic Galois group is as large as possible (i.e. with $\underline{\fu}(M)=   W_{-1}\inEnd(M)$). These motives are in particular interesting for the transcendence properties of their periods: in view of Grothendieck's period conjecture the field generated by their periods should have the highest possible transcendence degree among all motives with the same associated graded.
\medskip\par
A particularly striking implication of our result is that a suggestion of Euler about $\zeta(3)$ is incompatible with Grothendieck's period conjecture. 
In a 1785 paper \cite{Euler}, Euler speculated that there may be rational numbers $\alpha$ and $\beta $ and an expression of the form
$$
\zeta(3)\ =\ \alpha (\log 2)^3 \ +\ \beta \pi^2 (\log 2).
$$
See the article of Dunham \cite{Dunham} which gives a very readable account of this statement and Euler's remarkable work on evaluating the
Riemann zeta function at integer arguments. 
\medskip\par
In section  6.8, we construct a mixed Tate motive with periods (essentially) $\zeta(3)$, $\log 2$, $\pi$ and a fourth period. Moreover, we use our
results to show that the dimension of the Galois group in this case is $4$. Thus, the period conjecture would predict that these four periods are
algebraically independent, and this is incompatible with Euler's expectation stated above. A more detailed description of this mixed Tate motive is given below. 
\medskip\par
\subsection{$\underline{\fu}(M)$ and the extensions $\mathcal{E}_p(M)$}
To be more precise, $\underline{\fu}(M)$ is the subobject of $W_{-1}\inEnd(M)$ with the property that if $\omega$ is any fiber functor over $K$, then 
\[ \omega \, \underline{\fu}(M) \ \subset \ \omega \, W_{-1}\inEnd(M) \ = \ W_{-1} End(\omega M) \]
is the Lie algebra of
\[
\mathcal{U}(M,\omega) \ := \ \ker\bigm(\mathcal{G}(M,\omega) \ \xrightarrow{\text{restriction}} \ \mathcal{G}(Gr^WM,\omega)\bigm),
\]
where $\mathcal{G}(-,\omega)$ denotes the fundamental group of the indicated object with respect to $\omega$. If $Gr^WM$ is semisimple (which will be the case if $\mathbf{T}$ is a category of motives), then $\mathcal{U}(M,\omega)$ is the unipotent radical of $\mathcal{G}(M,\omega)$. 
\medskip\par
As stated above, Deligne (see \cite[Appendix]{Jo14}) describes $\underline{\fu}(M)$ in terms of extensions that arise naturally from the weight filtration on $M$: For each integer $p$, let $\mathcal{E}_p(M)$ be the {\it $p$-th extension class of $M$} given by Eq. \eqref{eq1into}, considered as an extension of the unit object $\mathbbm{1}$ by $\inHom(M/W_pM, W_pM)$. Pushing this extension forward along the natural injection
\[
\inHom(M/W_pM, W_pM) \ \longrightarrow \ W_{-1}\inEnd(M)
\] 
we get an extension of $\mathbbm{1}$ by $W_{-1}\inEnd(M)$, which we also denote by $\mathcal{E}_p(M)$. The {\it total} extension class of $M$ is then the extension 
\[
\mathcal{E}(M) \ := \ \sum\limits_p \mathcal{E}_p(M) \ \in \ Ext^1(\mathbbm{1}, W_{-1}\inEnd(M)).
\]
Deligne's result asserts that $\underline{\fu}(M)$ is the smallest subobject of $W_{-1}\inEnd(M)$ such that $\mathcal{E}(M)$ is in the image of the pushforward 
\begin{equation}\label{eq2intro}
Ext^1(\mathbbm{1}, \underline{\fu}(M)) \ \longrightarrow \ Ext^1(\mathbbm{1}, W_{-1}\inEnd(M))
\end{equation}
under the inclusion $\underline{\fu}(M)\subset W_{-1}\inEnd(M)$. Deligne proves this in part by exploiting the weight filtration to construct an explicit extension of $\mathbbm{1}$ by $\underline{\fu}(M)$ which pushes forward to $\mathcal{E}(M)$. 
\medskip\par
The first half of this paper is dedicated to the study of the relation between $\underline{\fu}(M)$ and the individual extensions $\mathcal{E}_p(M)
$, with a view to refining Deligne's result.  In general, the individual extensions $\mathcal{E}_p(M)$ may not be in the image of the pushforward 
map Eq. \eqref{eq2intro}; an example involving 1-motives can be given using the work \cite{JR86} of Jacquinot and Ribet on deficient points on 
semiabelian varieties (see Section \ref{counterexample to Thm3 without IA hypothesis} and the remarks at its end). The main result of the first 
half of the paper gives a sufficient condition for when the extension $\mathcal{E}_p(M)$ is in the image of Eq. \eqref{eq2intro} (see Theorem \ref{thm 2} and its corollaries). 
\medskip\par
\subsection{A more detailed overview}
We continue this introduction by giving a more detailed overview of the contents of the paper, starting with the first half. Fix an integer $p$ and an object $M$ of $\mathbf{T}$. It is natural to expect $\mathcal{E}_p(M)$ to be related to the subobject
\[
\underline{\fu}_p(M) \ := \ \underline{\fu}(M) \cap \inHom(M/W_pM,W_pM)
\]
of $\underline{\fu}(M)$, where we have considered $\inHom(M/W_pM,W_pM)$ as a subobject of $W_{-1}\inEnd(M)$ via the natural injection. This is indeed the case: Write $\mathcal{E}_p(M)$ explicitly as 
\begin{equation}\label{eq3intro}
0 \ \longrightarrow \ \inHom(M/W_pM, W_pM) \ \longrightarrow \ \inHom(M/W_pM, M)^\dagger \ \longrightarrow \ \mathbbm{1} \ \longrightarrow \ 0
\end{equation}
(see Section \ref{explicit description of E_p} for the explicit description of the middle object). Then by Theorem 3.3.1 of \cite{EM21} (which is proved by a small modification of Hardouin's proof of Theorem 2 of \cite{Har11}), we have\\
($\ast$) : { \it $\underline{\fu}_p(M)$ is the smallest subobject of $\inHom(M/W_pM,W_pM)$ such that 
\[
\inHom(M/W_pM, M)^\dagger \, / \underline{\fu}_p(M)
\]
belongs to the subcategory $\langle W_pM, M/W_pM\rangle^{\otimes}$.}\\
Here, as usual, the notation $\langle ~\rangle^{\otimes}$ means the smallest full Tannakian subcategory containing the indicated objects and closed under subobjects. The first contribution of the present article is to reformulate this statement in a more natural way, in the language of extensions originating from subcategories
(discussed in Section \ref{section extensions coming from subcategories}). Given a full Tannakian subcategory $\mathbf{S}$ of $\mathbf{T}$ which is closed under subobjects, we say an extension $\mathcal{E}$ of $\mathbbm{1}$ by an object $A$ of $\mathbf{T}$ originates from $\mathbf{S}$ if there is an object $A'$ of $\mathbf{S}$, an extension $\mathcal{E}'$ of $\mathbbm{1}$ by $A'$ in $\mathbf{S}$, and a morphism $A'\longrightarrow A$ under which $\mathcal{E}'$ pushes forward to  $\mathcal{E}$. While this is a very natural and simple generalization of the notion of splitting of sequences (as an extension splits if and only if it originates from a semisimple $\mathbf{S}$), it opens the door to refinements of Statement ($\ast$) and Deligne's theorem. The reformulated version of Statement ($\ast$) is given in Theorem \ref{thm 1}. It asserts that $\underline{\fu}_p(M)$ is the smallest subobject of $\inHom(M/W_pM,W_pM)$ such that the pushforward $\mathcal{E}_p(M)/\underline{\fu}_p(M)$ of $\mathcal{E}_p(M)$ under the quotient map originates from the subcategory 
\begin{equation}\label{eq106}
\langle W_pM, M/W_pM\rangle^{\otimes}.
\end{equation}
Note that one advantage of formulating the statement in this language is that here we may think of $\mathcal{E}_p(M)$ as an extension of $\mathbbm{1}$ by $\inHom(M/W_pM,W_pM)$ or by $W_{-1}\inEnd(M)$; see Section \ref{section:  part 1.1}.
\medskip\par
Our next goal is to find refinements of Theorem \ref{thm 1} in which the category Eq. \eqref{eq106} is replaced by smaller categories. Ideally, this category can be replaced by a semisimple category, in which case the pushforward $\mathcal{E}_p(M)/\underline{\fu}_p(M)$ of $\mathcal{E}_p(M)$ along the quotient map
will split. (By weight considerations and the long exact sequence for $Ext$ groups this is equivalent to $\mathcal{E}_p(M)$ being in the image of Eq. \eqref{eq2intro}.) But from the examples of 1-motives mentioned earlier we know that in general, this will not be the case.
\medskip\par
Let $q\leq p$. The second contribution of this paper is to show that if $M$ satisfies certain ``independence axioms", then in the statement of Theorem \ref{thm 1} the category Eq. \eqref{eq106} can be replaced by the smaller category $\langle W_qM, Gr^WM\rangle^{\otimes}$ (smaller because $q\leq p$); this is Theorem \ref{thm 2} in Section \ref{section: part 1.2}. The independence axioms are given in Section \ref{par: IA axioms}, and in fact, only depend on $Gr^WM$. Roughly speaking, they require the subobject
\begin{equation}\label{eq4intro}
\bigoplus_{\stackrel{i,j}{j>q, i}} \, \inHom(Gr^W_jM, Gr^W_iM)
\end{equation}
of $W_{-1}\inEnd(Gr^WM)$ to suitably decompose as a direct sum of two ``independent" summands. In the weak sense, here the word ``independent" means not having any nonzero isomorphic subobjects, and in the strong sense, it means having disjoint sets of weights (see the axioms $(IA1)_{\{p,q\}}$ and $(IA2)_{\{p,q\}}$ in Section \ref{par: IA axioms}). 
\medskip\par
An interesting consequence of Theorem \ref{thm 2} is the following refinement of Deligne's theorem (see Corollary \ref{cor 1 of thm 2}): if $Gr^WM$ is semisimple (e.g. if $\mathbf{T}$ is a category of motives) and the weak independence axioms hold for all $q\leq p$, then $\mathcal{E}_p(M)/\underline{\fu}_p(M)$ splits. In particular, if $Gr^WM$ is semisimple and $W_{-1}\inEnd(M)$ has ${n \choose 2}$ distinct weights where $n$ is the number of weights of $M$ (e.g. if $M$ has weights $0,-1,-3, -7$), then every $\mathcal{E}_p(M)/\underline{\fu}_p(M)$ splits for every $p$ (see Corollary \ref{cor 2 of thm 2}).
\medskip\par
The proof of Theorem \ref{thm 2} is similar to the proof of Statement ($\ast$) (or rather, of Theorem \ref{thm 1}), albeit with two added ingredients. Let $\underline{\fu}_{\geq q}(M)$ be the Lie algebra of the kernel of the restriction map from the fundamental group of $M$ to that of $W_qM\oplus Gr^WM$. The first new component is thanks to the independence axioms: they guarantee that $Gr^W \underline{\fu}_{\geq q}(M)$ (which is a subobject of Eq. \eqref{eq4intro}) decomposes according to the decomposition of Eq. \eqref{eq4intro} into our independent objects (see Lemma \ref{decomposition of Gr u_>=p}). This is the only place in the proof of Theorem \ref{thm 2} that the independence axioms play a part. Taking $\omega$ to be any fiber functor, this gives a decomposition of $\omega \, Gr^W \underline{\fu}_{\geq q}(M)$. The second added ingredient is that we use the fundamental theorem of Tannakian categories with $\omega\circ Gr^W$ as the fiber functor (rather than using $\omega$ itself). Notice the difference in the nature of this type of argument and Deligne's argument in \cite[Appendix]{Jo14}, which explicitly constructs an extension of $\mathbbm{1}$ by $\underline{\fu}(M)$ that pushes forward to the total class of $M$. We should point out that the idea of working with the associated graded fiber functor already appears in \cite{De94}, and since then has featured frequently in the literature, especially in the setting of categories of mixed Tate motives (e.g. \cite{DG05}).
\medskip\par
It would be interesting to give a more conceptual explanation (or geometric interpretation, in the case of motives) for the fact that the independence axioms force the individual extension classes $\mathcal{E}_p(M)/\underline{\fu}(M)$ (or $\mathcal{E}_p(M)/\underline{\fu}_p(M)$) to split. 
\medskip\par
We now discuss the contents of the second half of the paper (Section \ref{sec: part 2}). Let $\mathbf{T}$ be a Tannakian category of mixed motives over a field $K$ of characteristic zero, e.g. the Tannakian categories of Nori or Ayoub of mixed motives over $K$, or Voevodsky's Tannakian category of mixed Tate motives over a number field, or categories of mixed motives defined in terms of realizations. We say $\underline{\fu}(M)$ is {\it large} (or that $M$ has a large $\underline{\fu}$) if $\underline{\fu}(M)$ is equal to $W_{-1}\inEnd(M)$. As we pointed out earlier, such motives are interesting from the point of view of the transcendence properties of their periods.  Our original motivation for this part of the paper was to study (or ideally, classify up to isomorphism) all objects $M$ with large $\underline{\fu}$ and associated graded isomorphic to 
\[
\QQ(n) \oplus A \oplus \mathbbm{1},
\]
where $A$ is a given pure object of weight $p$ with $-2n<p<0$. We then realized that much of the discussion can be given in more generality, leading to the contents of this part of the paper as currently presented (and reviewed below). 
\medskip\par
Suppose tentatively that $M$ is an extension of $\mathbbm{1}$ by an object $L$ of highest weight $p$ with $p<0$. It is easy to see that if $\underline{\fu}(M)$ is large, then so are $\underline{\fu}(L)$ and $\underline{\fu}(M/W_{p-1}(L))$. The first main result of Section \ref{sec: part 2} (Theorem \ref{thm3}) gives a sufficient condition for the converse statement: it asserts that if $M$ satisfies a suitable independence axiom, and if $\underline{\fu}(L)$ and $\underline{\fu}(M/W_{p-1}(L))$ are large, then so is  $\underline{\fu}(M)$. This is an application of Corollary \ref{cor 1 of thm 2}. As in the case of the latter corollary, examples involving 1-motives show that the conclusion of Theorem \ref{thm3} is in general false without the hypothesis about the independence axiom (see Section \ref{counterexample to Thm3 without IA hypothesis}). 
\medskip\par
Theorem \ref{thm3} suggests a way to obtain more complicated objects with large $\underline{\fu}$ by ``patching together" smaller such objects. More precisely, given an object $L$ of highest weight $p$ with $p<0$ which has a large $\underline{\fu}$, and an object $N$ which is an extension of $\mathbbm{1}$ by $Gr^W_pL$ and also has a large $\underline{\fu}$, we can look for objects $M$ such that $W_pM\simeq L$ and $M/W_{p-1}M\simeq N$; assuming the relevant independence axiom (which only depends on $Gr^WM\simeq Gr^WL\oplus\mathbbm{1}$) holds, any such $M$ has a large $\underline{\fu}$. The answer to the question of existence of such $M$ is given by Grothendieck's formalism of {\it extensions panach\'{e}es} \cite{Gr68}: The obstruction is an element of $Ext^2(\mathbbm{1}, W_{p-1}L)$. Moreover, the object $M$ is unique up to isomorphism if $Ext^1(\mathbbm{1},W_{p-1}L)=0$ (see Lemma \ref{criterion for compatibility of a pair}).
\medskip\par
We consider the following classification problem in Sections \ref{sec: compatible pairs, defn} - \ref{sec: summary of results of the third part}: given $B$ of weights $< p$ and with a large $\underline{\fu}$, and a nonzero pure object $A$ of negative weight $p$, classify up to isomorphism all $M$ with large $\underline{\fu}$ satisfying $W_{p-1}M\simeq B$, $Gr_p^WM\simeq A$ and $M/W_pM\simeq \mathbbm{1}$ (with the isomorphisms not part of the data). We manage to give a complete solution to this problem when $B\oplus A\oplus \mathbbm{1}$ satisfies an independence axiom and $Ext^1(\mathbbm{1}, B)=0$; the solution is summarized in Section \ref{sec: summary of results of the third part}, just before Corollary \ref{cor: motives with three weights and Gr^W isomorphic to Q(n)+A+1, bijective case }. To get there, in Sections \ref{sec: compatible pairs, defn} - \ref{sec: classification result compatible pairs vs objects with ... up to isomorphism} we study the {\it extensions panach\'{e}es} problem\footnote{Actually, a slight variation of it; see the beginning of Section \ref{sec: compatible pairs, defn}.} in the setting of an abelian category with weights. The main result is summarized in Proposition \ref{prop compatible pairs} (see also Lemma \ref{equivalent pairs}). As a special case of these results, in Corollary \ref{cor: motives with three weights and Gr^W isomorphic to Q(n)+A+1, bijective case } we give an answer to our original motivating classification problem about objects with associated graded isomorphic to $\QQ(n) \oplus A \oplus \mathbbm{1}$. 
\medskip\par
In Section \ref{sec: classification of 3 dimensional mixed Tate motives with large u} we specialize to the category $\mathbf{MT}(\QQ)$ of (say, Voevodsky) mixed Tate motives over $\QQ$. The nice feature here is that the Ext groups are known. We use Corollary \ref{cor: motives with three weights and Gr^W isomorphic to Q(n)+A+1, bijective case } to give a complete classification, up to isomorphism, of all 3-dimensional mixed Tate motives over $\QQ$ with large $\underline{\fu}$ and associated graded isomorphic to $\QQ(n)\oplus \QQ(k) \oplus \mathbbm{1}$ with $n>k>0$ and $n\neq 2k$; the very last condition is the independence axiom in this situation. 
\medskip\par
Let us consider an example from Section \ref{sec: classification of 3 dimensional mixed Tate motives with large u} here. Let $r$ be an integer $>1$ and $N$ the Kummer 1-motive $[\ZZ\stackrel{1\mapsto r}{\longrightarrow} \mathbb{G}_m]$, considered as an object of $\mathbf{MT}(\QQ)$. Let $n$ be an even integer $\geq 4$, and $L$ an object which is a nontrivial extension of $\mathbbm{1}$ by $\QQ(n-1)$ (so with $(2\pi i)^{1-n}\zeta(n-1)$ as a period). Since $Ext^2$ groups vanish in $\mathbf{MT}(\QQ)$ and $Ext^1(\mathbbm{1},\QQ(n))=0$, the two objects  $L(1)$ and $N$ can be patched together to form an object\footnote{This object is denoted by $M_{n,r}$ in Section \ref{sec: classification of 3 dimensional mixed Tate motives with large u}.} $M$ of $\mathbf{MT}(\QQ)$, unique up to isomorphism, such that $W_{-2}M\simeq L(1)$ and $M/\QQ(n)\simeq N$. Moreover, $M$ satisfies the required independence axiom (as $n\neq 2$), so that it follows from Theorem \ref{thm3} that $\underline{\fu}(M)$ is large. According to Grothendieck's period conjecture, the field generated over $\QQ$ by the periods of $M$ should have transcendence degree equal to
\[
\dim (\mathcal{G}(M,\omega_B)) \ = \ \dim(\omega_B\, \underline{\fu}(M)) \, + \, \dim (\mathcal{G}(Gr^WM,\omega_B)) \ = \ 3 + 1 \ = \ 4
\]
(where $\omega_B$ = Betti realization). The nonzero entries of the period matrix of $M$ with respect to suitably chosen bases of de Rham and Betti realizations are $(2\pi i)^{-n}, (2\pi i)^{-n}  \zeta(n-1), (2\pi i)^{-1}$ (coming from $L(1)$), $(2\pi i)^{-1}\log(r), 1$ (coming from $N$), and a ``new period". So Grothendieck's period conjecture predicts that 
\[
2\pi i, \, \zeta(n-1), \, \log(r), \, \text{and the new period of $M$}
\]
must be algebraically independent over $\QQ$.
\medskip\par
The new period discussed above seems rather mysterious, and it would be very interesting to somehow\footnote{Ideally, one would like to do this by giving a geometric construction of $M$, but this may be too difficult especially when $n>4$. In general, giving geometric constructions of mixed Tate motives with a few weights is a difficult problem. See \cite[\S 1.4]{Br16}.} calculate it. When $r$ is 2 (or a power of it), $M$ is a mixed Tate motive over $\ZZ[1/2]$, and hence by Deligne's work \cite{De10} the new period will be a linear combination of alternating multiple zeta values, which one should be able to calculate using the formula of Goncharov and Brown (\cite{Go05} and \cite{Br12}) for the motivic coaction on iterated integrals.\footnote{This was told to us by Cl\'{e}ment Dupont.} On the other hand, for general $r$, at least a priori, the new period may not be an iterated integral on the projective line $\mathbb{P}^1$ minus $\{0,\infty\}\cup \mu_r$. (This is related to the question of whether the category of mixed Tate motives over $\ZZ[1/r]$ is generated by the fundamental groupoid of $\mathbb{P}^1\setminus(\{0,\infty\}\cup \mu_r)$, and for $r>2$ one expects the answer to this question to be in general negative. See Section 3 of \cite{DW16} for a discussion of this question.)
\medskip\par
After the discussion of 3-dimensional mixed Tate motives with large $\underline{\fu}$, in Section \ref{sec: 4 dimensional examples} we briefly consider some 4-dimensional examples; this leads again to some interesting questions about periods. One difference between the 4-dimensional and 3-dimensional examples is that in the former case (at least, a priori) one gets a family of motives with a large $\underline{\fu}$ when patching together a 3-dimensional $L$ and a 2-dimensional $N$. 
\medskip\par
We end this introduction with some words on the organization of the paper. In Section \ref{sec: prelim} we review some basic material about Tannakian categories. The notion of extensions originating from subcategories of a Tannakian category is discussed in Section \ref{sec: extensions originating from a subcategory}. Here we prove a few lemmas on this concept that will be useful throughout the paper. Starting from Section \ref{section:  part 1.1} we work in a Tannakian category with a weight filtration. In Section \ref{section:  part 1.1} we introduce the relevant objects and give the reformulation of Statement ($\ast$) (Theorem \ref{thm 1}). The goal of Section \ref{section: part 1.2} is to give the main results of the first part of the paper (Theorem \ref{thm 2} and its corollaries), in which we show that the independence axioms introduced in the same section result in refinements of Theorem \ref{thm 1} and Deligne's theorem. At the end of Section \ref{section: part 1.2} we also prove a variant of Theorem \ref{thm 2} for $q>p$ case (see Thoerem \ref{thm2'}). Section \ref{sec: part 2} contains the application to motives with large unipotent radicals of motivic Galois groups, as discussed above. We should point out that prior to Section \ref{sec: classification of 3 dimensional mixed Tate motives with large u} we use the term ``motive" only because we find it more suggestive: the discussion is valid in any Tannakian category with a weight filtration as long as the word ``motive" is interpreted as ``an object with a semisimple weight associated graded". In discussions where the Tate objects $\QQ(n)$ play a role, we also need to assume that there is a pure object $\QQ(1)$ of weight -2 such that the functor $-(1):=-\otimes \QQ(1)$ is invertible. Sections \ref{sec: classification of 3 dimensional mixed Tate motives with large u} and \ref{sec: 4 dimensional examples} take place in the setting of a Tannakian category of mixed Tate motives over $\QQ$ with the ``correct" Ext groups. Finally, Section \ref{counterexample to Thm3 without IA hypothesis} uses 1-motives to give counter-examples to several statements in the paper, if the hypotheses regarding the independence axioms are omitted.
\medskip\par
{\bf Acknowledgements.} We would like to thank Daniel Bertrand and Madhav Nori for a few insightful correspondences. We also thank Cl\'{e}ment Dupont for a very helpful correspondence about the motives $M_{n,r}$ of Section \ref{sec: classification of 3 dimensional mixed Tate motives with large u}, and for providing us with some valuable references. 

\section{Preliminaries on Tannakian categories}\label{sec: prelim}

The goal of this section is to review certain generalities about fundamental groups in Tannakian categories and fix some notation. None of the results in this section are new. The reader can refer to \cite{DM82} for the basics of Tannakian categories, for instance.

\subsection{Notation} 
The category of groups is denoted by $\mathbf{Groups}$. For any commutative ring $R$, we denote the category of $R$-modules (resp. commutative $R$-algebras) by $\mathbf{Mod}_R$ (resp. $\mathbf{Alg}_R$). We often denote the Hom and End groups in a category of modules simply by $Hom$ and $End$, with the coefficient ring being understood from the context.

Throughout, $K$ is a field of characteristic zero. If $V$ is a vector space over $K$, we denote the general linear group of $V$ by $GL(V)$; it is an algebraic group over $K$. If $G$ is an algebraic group over $K$, we denote the Lie algebra of $G$ by $Lie(G)$, and the category of finite-dimensional representation of $G$ (over $K$) by $\mathbf{Rep}(G)$. 

As usual, given a morphism $\alpha: \omega\longrightarrow \omega'$ of functors, for any object $M$ of the domain category the corresponding morphism $\omega M \longrightarrow \omega'M$ in the target category is denoted by $\alpha_M$.

Finally, in various contexts, we use the notation $f|_X$ for the restriction of $f$ to $X$ (whatever $f$ and $X$ are).

\subsection{} By a Tannakian category over $K$ we mean a neutral Tannakian category over $K$, i.e., in the language of \cite{DM82}, a rigid abelian $K$-linear tensor category with $K$ as the endomorphism algebra of the unit object, for which a fiber functor over $K$ ( = an exact faithful\footnote{Actually including faithfulness here is redundant, as it follows from the rest of the requirements. See  \cite[\S 2.10, 2.11]{De90}.} $K$-linear tensor functor from the category to $\mathbf{Mod}_K$) exists. 

If $\mathbf{T}$ is a Tannakian category over $K$ and $\omega: \mathbf{T}\longrightarrow \mathbf{Mod}_K$ 
is a fiber functor (over $K$), we denote the fundamental group of $\mathbf{T}$ with respect to $\omega$ by $\mathcal{G}(\mathbf{T},\omega)$ ( = $\underline{Aut}^{\otimes}(\omega)$ in the standard notation); thus (by the fundamental theorem of Tannakian categories) this is an affine group scheme over $K$ with
\begin{align*}
\mathcal{G}(\mathbf{T},\omega)(R) \ = \ & \text{the group of automorphisms of the functor}\\
&~\omega  \otimes 1_R: \mathbf{T}  \longrightarrow \mathbf{Mod}_K \longrightarrow \mathbf{Mod}_R\\
&~\text{respecting the tensor structures}
\end{align*}
for any $K$-algebra $R$. For any object $M$ of $\mathbf{T}$, we have a representation 
\[
\rho_M \ : \ \mathcal{G}(\mathbf{T},\omega) \ \longrightarrow \ GL(\omega M) \hspace{.2in} \sigma\mapsto \sigma_M
\]
and (agian by the fundamental theorem) the functor
\[
\mathbf{T} \ \longrightarrow \ \mathbf{Rep}( \mathcal{G}(\mathbf{T},\omega)) \hspace{.2in} M \mapsto (\omega M, \rho_M),
\]
which with abuse of notation we also denote by $\omega$, is an equivalence of categories. 

\subsection{} Let $\mathbf{T}$ be a Tannakian category over $K$ with unit object denoted by $\mathbbm{1}$. Let 
\[
\omega: \mathbf{T} \ \longrightarrow \ \mathbf{Mod}_K
\] 
be a fiber functor. For any full Tannakian subcategory $\mathbf{S}$ of $\mathbf{T}$ which is closed under taking subobjects (and hence subquotients), the inclusion $\mathbf{S}\subset \mathbf{T}$ gives a surjective restriction map
\[ \mathcal{G}(T,\omega) \ \longrightarrow \ \mathcal{G}(\mathbf{S},\omega |_{\mathbf{S}}) \]
(surjective because $\mathbf{S}$ is closed under taking subobjects, see \cite[Proposition 2.21]{DM82}). 

\subsection{} Given any objects $M_1,\ldots M_n$ of $\mathbf{T}$, let $\langle M_1,\ldots, M_n\rangle^{\otimes}$ be the Tannakian subcategory generated by $M_1,\ldots M_n$; by definition, $\langle M_1,\ldots, M_n\rangle^{\otimes}$ is the smallest full Tannakian subcategory of $\mathbf{T}$ which contains the $M_i$ and is closed under taking subobjects. Every object of $\langle M_1,\ldots, M_n\rangle^{\otimes}$ is obtained from $M_1,\ldots M_n$ and $\mathbbm{1}$ by finitely many iterations of taking direct sums, tensor products, duals, and subobjects. We have
\[
\langle M_1,\ldots, M_n\rangle^{\otimes} \ = \ \langle \bigoplus_{1\leq i\leq n} M_i\rangle^{\otimes}.
\]

\subsection{} Let $M$ be an object of $\mathbf{T}$. Given a fiber functor $\omega$ over $K$, we set 
\[
\mathcal{G}(M,\omega) \ := \ \mathcal{G}(\langle M \rangle^{\otimes}, \omega |_{\langle M \rangle^{\otimes}}) \ = \ \underline{Aut}^{\otimes}(\omega |_{\langle M \rangle^{\otimes}});
\]
we call this the fundamental group of $M$ with respect to $\omega$. Since every object of $\langle M\rangle^{\otimes}$ is obtained from $M$ and $\mathbbm{1}$ by finitely many iterations of taking direct sums, tensor products, duals and subobjects, the map
\[
\rho_M \ : \ \mathcal{G}(M,\omega) \ \longrightarrow \ GL(\omega M)
\] 
(sending $\sigma$ to $\sigma_M$) is injective. In particular, $\mathcal{G}(M,\omega)$ is an algebraic group over $K$. 

Let $\mathfrak{g}(M,\omega )$ be the Lie algebra of $\mathcal{G}(M,\omega)$. In view of the equivalence of categories
\[
\langle M \rangle^{\otimes} \ \longrightarrow \ \mathbf{Rep}(\mathcal{G}(M,\omega))
\]
given by $\omega$, the adjoint representation of $\mathcal{G}(M,\omega)$ defines an object $\underline{\fg}(M,\omega)$ in $\langle M\rangle ^{\otimes}$ such that 
\[
\omega \, \underline{\fg}(M,\omega) \ = \ \mathfrak{g}(M,\omega)
\]
as representations of $\mathcal{G}(M,\omega)$, where the $\mathcal{G}(M,\omega)$-action on $\omega \underline{\fg}(M,\omega)$ corresponds to $\underline{\fg}(M,\omega)$ (i.e. is $\rho_{\underline{\fg}(M,\omega)}$) and the $\mathcal{G}(M,\omega)$-action on $\mathfrak{g}(M,\omega)$ is given by the adjoint representation.

Identify $\mathcal{G}(M,\omega)$ as a subgroup of $GL(\omega M)$ via $\rho_M$. This identifies
\begin{equation}\label{eq1}
\fg(M,\omega) \ \subset \ Lie(GL(\omega M)) = End(\omega M).
\end{equation}
Denote $\inEnd(M):=\inHom(M,M)$ (the internal $Hom$ in $\mathbf{T}$). Then we can identify $\omega\inEnd(M) = End(\omega M)$, with the action of $\mathcal{G}(M,\omega)$ on $End(\omega M)$ corresponding to $\inEnd(M)$ being by conjugation. The inclusion Eq. \eqref{eq1} is compatible with the actions of $\mathcal{G}(M,\omega)$, making 
\[
\underline{\fg}(M,\omega) \ \subset \ \inEnd(M).
\]
%
%

\subsection{} For any object $N$ of $\langle M\rangle^{\otimes}$, let $\mathcal{G}(M, N,\omega)$ be the kernel of the surjection
\[
\mathcal{G}(M,\omega) \ \longrightarrow \ \mathcal{G}(N,\omega)
\]
induced by the inclusion $\langle N\rangle^{\otimes} \subset \langle M\rangle^{\otimes}$ (so for instance, $\mathcal{G}(M, \mathbbm{1},\omega)=\mathcal{G}(M,\omega)$). The Lie subalgebra  
\[
\fg(M,N,\omega) \ := \ Lie(\mathcal{G}(M, N,\omega))
\]
of $\fg(M,\omega)$ is invariant under the adjoint action of $\mathcal{G}(M,\omega)$, giving rise to a subobject
\[
\underline{\fg}(M,N,\omega) \ \subset \ \underline{\fg}(M,\omega) \ \subset \ \inEnd(M).
\]

\subsection{} The subobjects $\underline{\fg}(M,N, \omega)$ of $\inEnd(M)$ do not depend on the choice of the fiber functor $\omega$. More precisely, for every object $N$ of $\langle M\rangle^{\otimes}$, there is a canonical subobject 
\[\underline{\fg}(M,N) \ \subset \ \inEnd(M)\]
such that for every $\omega$ over $K$,
\[
\omega \, \underline{\fg}(M, N) \ = \ \fg(M,N,\omega) \ \subset \ End(\omega M). 
\]
This can be seen via the machinery of algebraic geometry over a Tannakian category (\cite[\S 5 ,\S 6 ]{De89}) and is well-known, but in the interest of keeping the paper more self-contained, here we include a proof:

\begin{prop}\label{prop 2 prelim}
Suppose $\omega$ and $\omega'$ are two fiber functors $\mathbf{T}\longrightarrow \mathbf{Mod}_K$. Then for any objects $M$ of $\mathbf{T}$ and $N$ of $\langle M\rangle^{\otimes}$,
\[
\underline{\fg}(M,N,\omega) \ = \ \underline{\fg}(M,N,\omega') 
\]
(as subobjects of $\inEnd(M)$).
\end{prop}

\begin{proof}
By a theorem of Deligne (see \cite[\S 1.12, 1.13]{De90}), there exists a $K$-algebra $R$ such that the two functors $\omega\otimes 1_R$ and $\omega'\otimes 1_R$ are isomorphic as $\otimes$-functors. Let 
\[
\alpha \ : \ \omega\otimes 1_R \ \longrightarrow \ \omega'\otimes 1_R
\]
be an isomorphism respecting the tensor structures. Then conjugation by $\alpha |_{\langle M\rangle^{\otimes}}$ gives an isomorphism
\[
c_\alpha \ : \ \mathcal{G}(M,\omega)_R \ \longrightarrow \ \mathcal{G}(M,\omega')_R.
\]
On the other hand, conjugation by 
\[
\alpha_M \ : \ \omega M\otimes 1_R \ \longrightarrow \ \omega' M\otimes 1_R
\]
gives an isomorphism
\[
c_{\alpha_M} \ : \ GL(\omega M)_R \ \longrightarrow \ GL(\omega' M)_R.
\]
The maps $c_\alpha$ and $c_{\alpha_M}$ are compatible with one another, i.e. we have a commutative diagram
\[
\begin{tikzcd}[remember picture]
   \mathcal{G}(M,\omega)_R  \arrow{r}{c_\alpha, \ \simeq}   & \mathcal{G}(M,\omega')_R \\
   GL(\omega M)_R \arrow{r}{c_{\alpha_M}, \ \simeq} & GL(\omega' M)_R ,
\end{tikzcd}
\begin{tikzpicture}[overlay,remember picture]
\path (\tikzcdmatrixname-1-1) to node[midway,sloped]{$\subset$}
(\tikzcdmatrixname-2-1);
\path (\tikzcdmatrixname-1-2) to node[midway,sloped]{$\subset$}
(\tikzcdmatrixname-2-2);
\end{tikzpicture}
\]
where the vertical inclusions are by the identifications via $\rho_M$ for $\omega$ and $\omega'$ (i.e. are given by $\sigma\mapsto \sigma_M$). Going to the Lie algebras by taking derivatives we get a commutative diagram
\[\begin{tikzcd}[remember picture, column sep=large]
\omega\underline{\fg}(M,\omega)\otimes R \ = \ \fg(M,\omega)\otimes R  \arrow{r}{Dc_\alpha, \ \simeq}   & \fg(M,\omega')\otimes R \ = \ \omega'\underline{\fg}(M,\omega')\otimes R \\
\omega\inEnd(M) \otimes R  = \ End(\omega M)\otimes R \arrow{r}{Dc_{\alpha_M}, \ \simeq} & End(\omega' M)\otimes R \ = \ \omega'\inEnd(M).
\end{tikzcd}
\begin{tikzpicture}[overlay,remember picture]
\path (\tikzcdmatrixname-1-1) to node[midway,sloped]{$\subset$}
(\tikzcdmatrixname-2-1);
\path (\tikzcdmatrixname-1-2) to node[midway,sloped]{$\subset$}
(\tikzcdmatrixname-2-2);
\end{tikzpicture}
\]
The horizontal arrow in the second row is again just conjugation by $\alpha_M$, so that
\[
Dc_{\alpha_M} \ = \ \alpha_{\inEnd(M)}.
\]
On recalling that $\underline{\fg}(M,\omega)$ is a subobject of $\inEnd(M)$ and by commutativity of the previous diagram, we get
\begin{equation}\label{eq12}
\omega'\underline{\fg}(M,\omega) \otimes R \ = \ \alpha_{\inEnd(M)}(\omega\underline{\fg}(M,\omega)\otimes R) \ = \  \omega'\underline{\fg}(M,\omega')\otimes R
\end{equation}
(as subspaces of $End(\omega' M)\otimes R$). This shows that 
\[\omega'\underline{\fg}(M,\omega) \ = \ \omega'\underline{\fg}(M,\omega')\]
and hence $\underline{\fg}(M,\omega)=\underline{\fg}(M,\omega')$.

If $N$ is any object of $\langle M\rangle^{\otimes}$, by considering the analogous map to $c_\alpha$ for $N$ one easily sees that $c_\alpha$ maps $\mathcal{G}(M,N, \omega)_R$ onto $\mathcal{G}(M,N, \omega')_R$.
Thus 
\[
Dc_\alpha (\fg(M,N, \omega)\otimes R) \ = \ \fg(M,N, \omega')\otimes R,
\]
and similarly to Eq. \eqref{eq12} we get 
\[\omega'\underline{\fg}(M,N,\omega) \ = \ \omega'\underline{\fg}(M,N,\omega')\]
as subspaces of $End(\omega'M)$.
\end{proof}

\section{Extensions originating from a subcategory}\label{section extensions coming from subcategories}\label{sec: extensions originating from a subcategory}

The goal of this section is to introduce and prove a few lemmas about the basic but useful notion of extensions originating from subcategories of Tannakian categories. This concept will provide a natural language for the results of the paper. As in the previous section, $K$ is a field of characteristic zero.

\subsection{} Let $G$ be an affine group scheme over $K$. Let $H$ be a subgroup of $G$. Let $V$ be an object of $\mathbf{Rep}(G)$. Denote by $V^H$ the ($K$-) subspace of $V$ which is fixed by $H$. More precisely,
\[
V^H \ := \ \{v\in V: \forall R\in \mathbf{Alg}_K,  \forall \sigma\in H(R), \ \sigma(v\otimes 1_R)=v\otimes 1_R\}.
\]
Suppose $H$ is normal in $G$. Then $V^H$ is a $G$-subrepresentation of $V$ (i.e. a subobject of $V$ in $\mathbf{Rep}(G)$).

Consider the functor 
\[
\mathbf{Rep}(G/H) \ \longrightarrow \ \mathbf{Rep}(G)
\]
which considers a representation of $G/H$ as a representation of $G$ via the morphism $G\longrightarrow G/H$. This functor identifies $\mathbf{Rep}(G/H)$ as the full subcategory of $\mathbf{Rep}(G)$ consisting of those representation of $G$ on which $H$ acts trivially. It is evident that for every object $V$ of $\mathbf{Rep}(G)$, the object $V^H$ is the largest subobject of $V$ which belongs to the subcategory $\mathbf{Rep}(G/H)$.

\subsection{}\label{section criterion for objects being in subcategories in terms of the action of the kernel of the restriction map} Let $\mathbf{T}$ be a Tannakian category over $K$, with $\omega$ a fiber functor $\mathbf{T}\longrightarrow \mathbf{Mod}_K$. Let $\mathbf{S}$ be a full Tannakian subcategory of $\mathbf{T}$ which is closed under taking subobjects. The inclusion $\mathbf{S}\subset \mathbf{T}$ gives a surjection
\begin{equation}\label{eq5}
\mathcal{G}(\mathbf{T},\omega) \ \longrightarrow \ \mathcal{G}(\mathbf{S}, \omega |_{\mathbf{S}}).
\end{equation}
Denote the kernel of this map by $\mathcal{H}$.

Using the map Eq. \eqref{eq5} we may identify the category $\mathbf{Rep}(\mathcal{G}(\mathbf{S}, \omega |_{\mathbf{S}}))$ as the full subcategory of $\mathbf{Rep}(\mathcal{G}(\mathbf{T},\omega))$ consisting of all the objects on which $\mathcal{H}$ acts trivially. One has a commutative diagram 
\begin{equation}\label{eq11}
\begin{tikzcd}[remember picture]
   \mathbf{S} \arrow{r}{\omega |_{\mathbf{S}}, \ \simeq}   & \mathbf{Rep}(\mathcal{G}(\mathbf{S}, \omega |_{\mathbf{S}})) \\
    \mathbf{T}  \arrow{r}{\omega , \ \simeq} & \mathbf{Rep}(\mathcal{G}(\mathbf{T},\omega))  ,
\end{tikzcd}
\begin{tikzpicture}[overlay,remember picture]
\path (\tikzcdmatrixname-1-1) to node[midway,sloped]{$\subset$}
(\tikzcdmatrixname-2-1);
\path (\tikzcdmatrixname-1-2) to node[midway,sloped]{$\subset$}
(\tikzcdmatrixname-2-2);
\end{tikzpicture}
\end{equation}
where the horizontal arrows are the equivalences of categories given by the fundamental theorem of Tannakian categories. On recalling that $\mathbf{S}$ is closed under subobjects and hence in particular isomorphisms, it follows that any object $A$ of $\mathbf{T}$ belongs to the subcategory $\mathbf{S}$ if and only if $\mathcal{H}$ acts trivially on $\omega A$.

\subsection{} Let $A$ be an object of $\mathbf{T}$. Then $(\omega A)^\mathcal{H}$ is a $\mathcal{G}(\mathbf{T},\omega)$-subrepresentation of $\omega A$; hence there is a canonical subobject
\[
A_{\mathbf{S}} \ \subset \ A
\]
such that 
\[\omega(A_{\mathbf{S}})=(\omega A)^\mathcal{H}.\]
Since $(\omega A)^{\mathcal{H}}$ is the largest subobject of $\omega A \in \mathbf{Rep}(\mathcal{G}(\mathbf{T},\omega))$ which belongs to the subcategory $\mathbf{Rep}(\mathcal{G}(\mathbf{S}, \omega |_{\mathbf{S}}))$, it follows that $A_{\mathbf{S}}$ is the largest subobject of $A$ which belongs to $\mathbf{S}$.

Taking $\mathcal{H}$-invariance gives a left exact functor 
\[
\mathbf{Rep}(\mathcal{G}(\mathbf{T},\omega)) \ \longrightarrow \ \mathbf{Rep}(\mathcal{G}(\mathbf{S}, \omega |_{\mathbf{S}})).
\]
Thus we have a left exact functor 
\[
-_{\mathbf{S}} \ : \ \mathbf{T} \ \longrightarrow \ \mathbf{S}
\]
which on objects acts like $A\mapsto A_{\mathbf{S}}$ (and on morphisms acts by restriction of domain and codomain).

\subsection{}\label{section definition of extensions coming from subcategories and lemmas} Let $A$ be an object of $\mathbf{T}$. Let $\mathcal{E}$ in $Ext^1_\mathbf{T}(\mathbbm{1}, A)$ ( = Yoneda $Ext^1$ group in $\mathbf{T}$) be the class of the short exact sequence
\[
0 \ \longrightarrow \ A \ \longrightarrow \ E  \ \longrightarrow \ \mathbbm{1}  \ \longrightarrow \ 0.
\]
We say the extension $\mathcal{E}$  {\it originates from} or {\it comes from} $\mathbf{S}$ if there is a commutative diagram in $\mathbf{T}$
\begin{equation}\label{eq6}
\begin{tikzcd}
   0 \arrow[r] & A'  \arrow[d, ] \arrow[r, ] & E' \arrow[d, ] \arrow[r, ] &  \mathbbm{1} \ar[equal]{d} \arrow[r] & 0 \\
   0 \arrow[r] & A \arrow[r, ] & E \arrow[r, ] &  \mathbbm{1}  \arrow[r] & 0 ,
\end{tikzcd}
\end{equation}
where the rows are exact and the objects in the top row are in $\mathbf{S}$. In other words, we say $\mathcal{E}$ originates from $\mathbf{S}$ if there is an object $A'$ of $\mathbf{S}$ and a morphism $A'\longrightarrow A$ such that $\mathcal{E}$ is in the image of the pushforward map
\[
Ext^1_\mathbf{S}(\mathbbm{1}, A') \ \longrightarrow \ Ext^1_\mathbf{T}(\mathbbm{1},A).
\]

We now give a few lemmas on the notion of extensions originating from subcategories which are useful in the later sections. The lemmas take place in the above setting (i.e. with $\mathcal{E}$, $\mathbf{S}$, and $\mathcal{H}$ as above). The first lemma highlights that the notion of extensions originating from subcategories is a generalization of the notion of splitting of sequences.

\begin{lemma}
The following statements are equivalent:
\begin{itemize}
\item[(i)] The extension $\mathcal{E}$ splits.
\item[(ii)] The extension $\mathcal{E}$ originates from some semisimple $\mathbf{S}$.
\item[(iii)] The extension $\mathcal{E}$ originates from every $\mathbf{S}$.
\end{itemize}
\end{lemma}

\begin{proof}
The implications $(iii)\Longrightarrow (ii) \Longrightarrow (i)$ are trivial. As for $(i)\Longrightarrow (iii)$, note that if $\mathcal{E}$ splits, then it is the pushforward of the extension 
\[
0 \ \longrightarrow \ 0  \ \longrightarrow \ \mathbbm{1}  \ \longrightarrow \ \mathbbm{1}  \ \longrightarrow \ 0.
\]
\end{proof}

\begin{lemma}\label{lem 1 extensions originating from subcategories}
The following statements are equivalent:
\begin{itemize}
\item[(i)] The extension $\mathcal{E}$ originates from $\mathbf{S}$.
\item[(ii)] The extension $\omega\mathcal{E}$
\[
0 \ \longrightarrow \ \omega A \ \longrightarrow \ \omega E  \ \longrightarrow \ K  \ \longrightarrow \ 0
\]
splits in the category of representations of $\mathcal{H}$. 
\item[(iii)] The sequence
\[
0 \ \longrightarrow \ (\omega A)^\mathcal{H} \ \longrightarrow \ (\omega E)^\mathcal{H}  \ \longrightarrow \ K  \ \longrightarrow \ 0
\] 
(obtained by applying $\mathcal{H}$-invariance to $\omega \mathcal{E}$) is exact. 
\item[(iv)] The sequence in $\mathbf{S}$
\[
0 \ \longrightarrow \ A_\mathbf{S} \ \longrightarrow \ E_\mathbf{S}  \ \longrightarrow \ \mathbbm{1}  \ \longrightarrow \ 0
\]
obtained by applying $-_\mathbf{S}$ to the defining sequence of $\mathcal{E}$ is exact.
\end{itemize}
\end{lemma}

\begin{proof}

The equivalence of (iii) and (iv) is clear, as the sequence in (iii) is obtained by applying $\omega$ to the sequence in (iv). Note that since the functors $-^\mathcal{H}$ and $-_\mathbf{S}$ are left exact, the statements in (iii) and (iv) are really just statements about surjectivity of $(\omega E)^\mathcal{H}\longrightarrow K$ and $E_\mathbf{S}\longrightarrow \mathbbm{1}$. The implication (iv) $\Longrightarrow$ (i) is also clear, as we can use the extension given in (iv) as the top row in Eq. \eqref{eq6}. 

(i) $\Longrightarrow$ (iv): Suppose $\mathcal{E}$ originates from $\mathbf{S}$, with a commutative diagram as in Eq. \eqref{eq6}, with exact rows and the top row in $\mathbf{S}$. Since $\mathbf{S}$ is closed under taking subquotients, by replacing $A'$ and $E'$ if necessary by their images in $A$ and $E$, we may assume without loss of generality that $A'\subset A$ and $E'\subset E$, with the vertical arrows being considered as inclusion maps. Since $E'$ is in $\mathbf{S}$, we have $E'\subset E_\mathbf{S}$. This proves that the restriction of the surjection $E\longrightarrow \mathbbm{1}$ to $E_\mathbf{S}$ is still surjective, thus giving (iv).

(iii) $\Longrightarrow$ (ii): There is a commutative diagram of $\mathcal{G}(\mathbf{T},\omega)$-representations
\[
\begin{tikzcd}
   0 \arrow[r] & (\omega A)^\mathcal{H} \arrow[d, ] \arrow[r, ] & (\omega E')^\mathcal{H} \arrow[d, ] \arrow[r, ] &  K \ar[equal]{d} \arrow[r] & 0 \\
   0 \arrow[r] & \omega A \arrow[r, ] & \omega E \arrow[r, ] &  K  \arrow[r] & 0 ,
\end{tikzcd}
\]
where the bottom row is $\omega\mathcal{E}$, the vertical arrows are inclusion, and the rows are exact. Consider this diagram in the category of representations of $\mathcal{H}$. The top row splits, hence so does the bottom row ( = the pushout of the top row).

(ii) $\Longrightarrow$ (iii): Suppose (ii) holds. Choose a section $s$ of $\omega E\longrightarrow K$ in $\mathbf{Rep}(\mathcal{H})$. Then $s(1)$ is fixed by $\mathcal{H}$ and thus belongs to $(\omega E)^\mathcal{H}$. It follows that $(\omega E)^\mathcal{H}\longrightarrow K$ is surjective. 

\end{proof}

\begin{lemma}\label{lem 2 extensions originating from subcategories}
Suppose $A$ is an object of $\mathbf{S}$. Then $\mathcal{E}$ originates from $\mathbf{S}$ if and only if $E$ is an object of $\mathbf{S}$.
\end{lemma}

\begin{proof}
The ``if" implication is trivial. As for the ``only if" implication, suppose we have a diagram as in Eq. \eqref{eq6}, with exact rows and the objects of the top row in $\mathbf{S}$. Then $E$ is isomorphic to the fibered coproduct of $A$ and $E'$ over $A'$. Since $A$ and $E'$ are in $\mathbf{S}$, so is $E$.
\end{proof}

\begin{lemma}\label{lem 3 extensions originating from subcategories}
Let $A'$ be a subobject of $A$ such that the pushforward map
\[
Ext^1_{\mathbf{T}}(\mathbbm{1}, A'+A_{\mathbf{S}}) \ \longrightarrow \ Ext^1_{\mathbf{T}}(\mathbbm{1}, A)
\]
(along the inclusion $A'+A_{\mathbf{S}}\longrightarrow A$) is injective. Suppose $\mathcal{E}$ is the pushforward of an extension 
\[
\mathcal{E}' \ \in \ Ext^1_{\mathbf{T}}(\mathbbm{1}, A') 
\]
along the inclusion map $A'\longrightarrow A$. Then $\mathcal{E}$ originates from $\mathbf{S}$ if and only if $\mathcal{E}'$ does.
\end{lemma}

\begin{proof}
If $\mathcal{E}'$ originates from $\mathbf{S}$, then clearly so does $\mathcal{E}$. Suppose $\mathcal{E}$ originates from $\mathbf{S}$. Then $\mathcal{E}$ is the pushforward of the extension $\mathcal{E}_{\mathbf{S}}$ given in Statement (iv) of Lemma \ref{lem 1 extensions originating from subcategories} under the inclusion $A_{\mathbf{S}}\longrightarrow A$. Let $i: A_{\mathbf{S}}\longrightarrow A'+A_{\mathbf{S}}$ and $i': A'\longrightarrow A'+A_{\mathbf{S}}$ be inclusion maps. Apply the $\delta$-functor $Hom_{\mathbf{T}}(\mathbbm{1}, -)$ to the short exact sequence
\[
0 \ \longrightarrow \ A'\cap A_{\mathbf{S}} \ \longrightarrow \ A'\oplus A_{\mathbf{S}}   \ \stackrel{i-i'}{\longrightarrow} \ A' + A_{\mathbf{S}} \ \longrightarrow \ 0
\]
(where the injective arrow is the diagonal embedding). We get exact 
\[
Ext^1_{\mathbf{T}}(\mathbbm{1}, A'\cap A_{\mathbf{S}}) \ \longrightarrow \ Ext^1_{\mathbf{T}}(\mathbbm{1}, A')\oplus Ext^1_{\mathbf{T}}(\mathbbm{1}, A_{\mathbf{S}}) \ \stackrel{i_\ast-i'_\ast}{\longrightarrow} \ Ext^1_{\mathbf{T}}(\mathbbm{1}, A' + A_{\mathbf{S}}),
\]
where the lower stars denote pushforwards. The pushforward of the extension
\[
i_\ast(\mathcal{E}_{\mathbf{S}})-i'_\ast(\mathcal{E}') \ \in \ Ext^1_{\mathbf{T}}(\mathbbm{1}, A' + A_{\mathbf{S}})
\]
in $Ext^1_{\mathbf{T}}(\mathbbm{1}, A)$ is zero. By the injectivity hypothesis in the statement, $i_\ast(\mathcal{E}_S)-i'_\ast(\mathcal{E}')$ is already zero. It follows that there is an extension $\mathcal{E}''$
\[
0 \ \longrightarrow \ A'\cap A_{\mathbf{S}} \ \longrightarrow \ E''   \ \longrightarrow \ \mathbbm{1} \ \longrightarrow \ 0
\]
which pushes forward (under inclusion maps) to both $\mathcal{E}'$ and $\mathcal{E}_{\mathbf{S}}$. But then $A'\cap A_{\mathbf{S}}$ and $E''$, being subobjects of $A_{\mathbf{S}}$ and $E_{\mathbf{S}}$, belong to $\mathbf{S}$. Since $\mathcal{E}''$ pushes forward to $\mathcal{E}'$, the latter extension originates from $\mathbf{S}$.
\end{proof}

\begin{rem}
Note that the injectivity hypothesis in the statement of the previous lemma is guaranteed if 
\[
Hom_{\mathbf{T}}(\mathbbm{1}, A/(A'+A_{\mathbf{S}})) \ = \ 0
\] 
(and this will be the case whenever we use the result in the paper). This can be seen from the long exact sequence obtained by applying $Hom_{\mathbf{T}}(\mathbbm{1}, -)$ to 
\[
0 \ \longrightarrow \ A'+A_{\mathbf{S}} \ \longrightarrow \ A   \ \longrightarrow \ A/(A' + A_{\mathbf{S}}) \ \longrightarrow \ 0.
\]
\end{rem}

\section{Extension classes and subgroups of the fundamental group - Part I}\label{section:  part 1.1}

\subsection{}\label{par: part 1 setup} From this point on we suppose that $\mathbf{T}$ is a Tannakian category over a field $K$ of characteristic zero, equipped with a functorial exact finite increasing filtration $W_\cdot$, compatible with the tensor structure. We refer to $W_\cdot$ as the weight filtration. Here, the expression ``functorial exact finite increasing filtration $W_\cdot$" means that for every integer $n$, we have an exact functor $W_n: \mathbf{T} \longrightarrow \mathbf{T}$, such that for every object $M$ of $\mathbf{T}$, we have 
\begin{align*}
W_{n-1}M \ &\subset \ W_nM \hspace{.3in}(\forall n)\\
W_nM \ &= \ 0 \hspace{.3in}(\forall n\ll0)\\
W_nM \ &= \ M \hspace{.3in}(\forall n\gg 0),
\end{align*}
and such that the inclusions $W_nM\subset M$ for various $M$ give a morphism of functors from $W_n$ to the identity (and hence the $W_n$ form an inductive system of functors). Compatibility with the tensor product means that for every objects $M$ and $N$, we have
\begin{equation}\label{eq103}
W_n (M\otimes N) \ = \ \sum\limits_{\stackrel{p,q}{p+q=n}}\, W_pM\otimes W_qN.
\end{equation}
The associated graded functor $Gr^W$ is the functor defined on objects by 
\[Gr^WM \ := \ \bigoplus_n Gr^W_nM,\] 
where $Gr^W_nM:=W_nM/W_{n-1}M$, and on morphisms in the obvious way using the fact that we have morphisms of functors $W_{n-1}\longrightarrow W_n$. By the snake lemma, the associated graded functor (in fact, each $Gr^W_n$) is also exact. Also $Gr^W$ is a graded tensor functor, in the sense that (via a canonical isomorphism) we have
\[
Gr^W(M\otimes N) \ = \ Gr^W(M)\otimes Gr^W(N),
\]
with this identification being compatible with weights, i.e. being the direct sum of identifications
\[
Gr^W_n (M\otimes N) \ = \ \bigoplus\limits_{\stackrel{p,q}{p+q=n}}\, Gr^W_pM\otimes Gr^W_qM
\]
induced by Eq. \eqref{eq103}.

As it is customary, we call an object $M$ with $W_{n-1}M=0$ and $W_nM=M$ a pure object of weight $n$. Note that unless otherwise indicated, we do not assume that an object of the form $Gr^WM$ (i.e. a direct sum of pure objects) is necessarily semisimple.

Given any fiber functor $\omega$ (over $K$) and any object $M$, set
\[
W_\cdot \omega M \ := \ \omega(W_\cdot M).
\]
This defines an exact $\otimes$-filtration on $\omega$, in the language of Saavedra Rivano \cite{Riv72}, Chapter IV, \S 2 (note that Saavedra Rivano works with decreasing filtrations instead, and that his Condition FE 1) is guaranteed here because $K$ is a field). 

Given any objects $M$ and $N$, we identify 
\[
\omega \, \inHom(M,N) \ = \ Hom(\omega M, \omega N).
\]
One can then show that
\[
\omega \, W_n\inHom(M,N) \ = \ \{f\in Hom(\omega M, \omega N)\, : \, f( W_\cdot \omega M)\subset W_{\cdot+n}\omega N\}.
\]

\subsection{} Here and elsewhere in the paper, we shall use the notation and conventions of Section \ref{sec: prelim} for Tannakian fundamental groups and their Lie algebras.

Let $M$ be an object of $\mathbf{T}$. Given any fiber functor $\omega$, let $G(M,\omega)$ be the parabolic subgroup of $GL(\omega M)$ which stabilizes the filtration $W_\cdot$. Then
\[
Lie(G(M,\omega)) \ = \ W_0 End(\omega M).
\]
The elements of $\mathcal{G}(M,\omega)$ ( = the fundamental group of $M$ with respect to $\omega$) preserve subobjects of $M$, so that 
\[
\mathcal{G}(M,\omega) \ \subset \ G(M,\omega).
\]  
Going to the Lie algebras we have
\[
\underline{\fg}(M) \ \subset \ W_0\inEnd(M) . 
\] 
Every element of $G(M,\omega)$ induces an automorphism of $Gr^W\omega M$, giving rise to a homomorphism
\[
G(M,\omega) \ \longrightarrow \ GL(Gr^W\omega M).
\] 
Let $U(M,\omega)$ be the kernel of this map; then $U(M,\omega)$ is the unipotent radical of $G(M,\omega)$. It is easy to see that
\[
Lie(U(M,\omega)) \ = \ W_{-1} End(\omega M).
\]
Set
\[
\mathcal{U}(M,\omega) \ := \ \mathcal{G}(M, Gr^WM,\omega)
\]
( = the kernel of the restriction map $\mathcal{G}(M,\omega) \ \longrightarrow \ \mathcal{G}(Gr^WM,\omega)$ induced by the inclusion $\langle Gr^WM\rangle^{\otimes}\subset \langle M\rangle^{\otimes}$). Then
\begin{equation}\label{eq100}
\mathcal{U}(M,\omega) \ = \ \mathcal{G}(M,\omega)\cap U(M,\omega).
\end{equation}
In particular, $\mathcal{U}(M,\omega)$ is a unipotent group. If $\mathcal{G}(Gr^WM,\omega)$ happens to be reductive (i.e. if $Gr^WM$ is semisimple), then $\mathcal{U}(M,\omega)$ will be the unipotent radical of $\mathcal{G}(M,\omega)$.

We set
\[
\underline{\fu}(M) \ := \ \underline{\fg}(M, Gr^WM)
\]
and 
\[
\fu(M,\omega) \ := \ Lie \, \mathcal{U}(M,\omega)
\]
( = $ \fg(M, Gr^WM, \omega)$ in the notation of Section \ref{sec: prelim}). Then (for every $\omega$),
\[
\omega\, \underline{\fu}(M) \ = \ \fu(M,\omega).
\]

By Eq. \eqref{eq100}, we have
\[
\underline{\fu}(M) \ = \ \underline{\fg}(M)) \cap W_{-1}\inEnd(M)
\]
(where the intersection means the fibered product of the inclusions in $\inEnd(M)$).

\subsection{}\label{Deligne review} A result of Deligne (written by Jossen in the appendix of \cite{Jo14}) describes the subobject $\underline{\fu}(M)$ of $W_{-1}\inEnd(M)$ as follows\footnote{We thank Peter Jossen for patiently explaining to us some parts of Deligne's argument from \cite[Appendix]{Jo14}.}. From now on, if there is no ambiguity, we shall simply write $Hom$ (resp. $Ext^i$) for the Hom groups $Hom_\mathbf{T}$ (resp. the Yoneda $Ext^i_\mathbf{T}$ groups) in $\mathbf{T}$.

Recall from the Introduction that for each integer $p$, the $p$-th extension class 
\[
\mathcal{E}_p(M) \ \in \ Ext^1(\mathbbm{1},\inHom(M/W_pM, W_pM))
\]
of $M$ is the extension corresponding to the sequence
\[
0 \ \longrightarrow \ W_pM \ \longrightarrow \ M \ \longrightarrow \ M/W_pM \ \longrightarrow \ 0 
\]
under the canonical isomorphism
\begin{equation}\label{eq4}
Ext^1(M/W_pM, W_pM) \ \cong \ Ext^1(\mathbbm{1},\inHom(M/W_pM, W_pM)).
\end{equation}
Applying $\inHom(M/W_pM, -)$ to the inclusion $W_pM \longrightarrow M$ we get an injection
\[
\inHom(M/W_pM, W_pM) \ \longrightarrow \  \inHom(M/W_pM, M).
\]
On the other hand, applying $\inHom(-, M)$ to the quotient map $M\longrightarrow M/W_pM$ we get an injection
\[
\inHom(M/W_pM, M) \ \longrightarrow \ \inEnd(M).
\]
Composing the two injections, we get a map
\begin{equation}\label{eq2}
\inHom(M/W_pM, W_pM) \ \longrightarrow \ \inEnd(M).
\end{equation}
After applying a fiber functor $\omega$, this simply sends an element 
\[f\in Hom(\omega M/\omega W_pM, \omega W_pM)\] 
to the composition
\begin{equation}\label{eq7}
\omega M \ \xrightarrow{\text{ quotient }} \ \omega M/\omega W_pM \ \stackrel{f}{\longrightarrow} \ \omega W_pM \ \xrightarrow{\text{ inclusion }} \ \omega M.
\end{equation}
From this it is clear that indeed, the image of the map Eq. \eqref{eq2} is contained in $W_{-1}\inEnd(M)$. We shall identify $\inHom(M/W_pM, W_pM)$ as a subobject of $W_{-1}\inEnd(M)$ via the map Eq. \eqref{eq2}. Note that $\inHom(M/W_pM, W_pM)$ is an abelian Lie subalgebra of $W_{-1}\inEnd(M)$.

Pushing forward extensions along the inclusion map we get a map
\begin{equation}\label{eq3}
Ext^1(\mathbbm{1}, \inHom(M/W_pM, W_pM)) \ \longrightarrow \ Ext^1(\mathbbm{1}, W_{-1}\inEnd(M)),
\end{equation}
which is injective, as (by weight considerations), 
\[
Hom(\mathbbm{1}, \frac{W_{-1}\inEnd(M)}{\inHom(M/W_pM, W_pM)}) \ = \ 0.
\]
To simplify the notation, we shall identify 
\[
Ext^1(\mathbbm{1}, \inHom(M/W_pM, W_pM))
\]
with its image under Eq. \eqref{eq3}. 

Deligne defines the (total) extension class of $M$ to be
\[
\mathcal{E}(M) \ := \ \sum\limits_p \mathcal{E}_p(M) \ \in \ Ext^1(\mathbbm{1}, W_{-1}\inEnd(M)),
\]
(this is denoted by $cl(M)$ in \cite{Jo14}), and proves that the extension $\mathcal{E}(M)$ can be used to describe $\underline{\fu}(M)$:

\begin{thm}[Deligne, Appendix of \cite{Jo14}]\label{Deligne thm} \ \\
$\underline{\fu}(M)$ is the smallest subobject of $W_{-1}\inEnd(M)$ such that the extension $\mathcal{E}(M)$ is the pushforward of an element of $Ext^1(\mathbbm{1}, \underline{\fu}(M))$ under the inclusion $\underline{\fu}(M)\longrightarrow W_{-1}\inEnd(M)$.
\end{thm}

By weight considerations, the pushforward map 
\begin{equation}\label{eq101}
Ext^1(\mathbbm{1}, \underline{\fu}(M)) \ \longrightarrow \ Ext^1(\mathbbm{1}, W_{-1}\inEnd(M))
\end{equation}
is injective, so that the element pushing forward to $\mathcal{E}(M)$ is indeed unique. 

\begin{rem}
As we pointed out in the Introduction, in general, the individual extensions $\mathcal{E}_p(M)$ may not be in the image of the pushforward map Eq. \eqref{eq101}. See Section \ref{counterexample to Thm3 without IA hypothesis} (and Remark (2) therein) for examples in the category of mixed Hodge structures using the Jacquinot-Ribet defficient points on semiabelian varieties. 
\end{rem}

\subsection{} We adopt the following notation for pushforwards of extensions along quotient maps. If $\mathcal{E}$ is an extension of an object $A$ by $B$, then for any subobject $B'$ of $B$ we denote the pushforward of $\mathcal{E}$ along the quotient $B\longrightarrow B/B'$ by $\mathcal{E}/B'$.

Given any subobject $A\subset W_{-1}\inEnd(M)$, applying the functor $Hom(\mathbbm{1}, -)$ to the short exact sequence 
\[
0 \ \longrightarrow \ A  \ \longrightarrow \  W_{-1}\inEnd(M) \ \longrightarrow \  W_{-1}\inEnd(M)/A \ \longrightarrow \ 0
\]
we get a long exact sequence. In particular, we have exact
\[
Ext^1(\mathbbm{1}, A) \ \longrightarrow \ Ext^1(\mathbbm{1}, W_{-1}\inEnd(M))  \ \longrightarrow \ Ext^1(\mathbbm{1}, W_{-1}\inEnd(M)/A), 
\]
where the arrows are pushforwards along inclusion and quotient maps. Thus Deligne's result can be equivalently stated as that $\underline{\fu}(M)$ is the smallest subobject of $W_{-1}\inEnd(M)$ such that the pushforward 
\[
\mathcal{E}(M) /\underline{\fu}(M) \ \in  \ Ext^1(\mathbbm{1}, W_{-1}\inEnd(M)/\underline{\fu}(M))
\]
of $\mathcal{E}(M)$ splits. 

The formulation of Theorem \ref{Deligne thm} as given in the statement is more natural for Deligne's proof, as his argument goes by constructing an explicit extension of $\mathbbm{1}$ by $\underline{\fu}(M)$ which pushes forward to $\mathcal{E}(M)$. The formulation in terms of $\mathcal{E}(M) /\underline{\fu}(M)$ is however more natural when one wants to study the individual extensions $\mathcal{E}_p(M)$, as we shall see.

\subsection{}\label{explicit description of E_p} The canonical isomorphism Eq. \eqref{eq4} is given by first applying the functor $\inHom(M/W_pM, -)$ to an element of $Ext^1(M/W_pM, W_pM)$, and then pulling back along the canonical map
\[
\mathbbm{1} \ \longrightarrow \ \inEnd(M/W_pM)
\]
(which after applying a fiber functor $\omega$, sends $1$ to the identity map on $\omega(M/W_pM)$). Going through this, we see that assuming $M/W_pM \neq 0$, the extension 
\[\mathcal{E}_p \ \in Ext^1(\mathbbm{1},  \inHom(M/W_pM, W_pM))\]
is the class of 
\begin{equation}\label{eq8}
0 \ \longrightarrow \ \inHom(M/W_pM, W_pM) \ \longrightarrow \ \inHom(M/W_pM, M)^\dagger  \ \longrightarrow \ \mathbbm{1}  \ \longrightarrow \ 0,
\end{equation}
where $\inHom(M/W_pM,M)^\dagger$ is the subobject of $\inHom(M/W_pM,M)$ characterized by
\begin{align*}
\omega \inHom(M/W_pM,M)^\dagger \ &= \ Hom(\omega M/\omega W_pM,\omega M)^\dagger \\ 
&:= \ \bigm\{f\in Hom(\omega M/\omega W_pM ,\omega M) : \\
& \ \ \ \ \ \ \ f\pmod{\omega W_pM} = \lambda(f) Id_{\omega M/\omega W_pM} \\ 
& \ \ \ \ \ \ \ \text{for some} \ \lambda(f)\in K\bigm\}
\end{align*}
for any fiber functor $\omega$. The injective (resp. surjective) arrow in Eq. \eqref{eq8} is, after applying $\omega$, the natural inclusion (resp. the map $f\mapsto \lambda(f)$ , with $\lambda(f)\in K$ as in the definition of $\inHom(M/W_pM,M)^\dagger$ above).

If $M/W_pM = 0$, then for convenience we set $\inHom(M/W_pM,M)^\dagger:=\mathbbm{1}$, so that $\mathcal{E}_p$ is still given by the sequence Eq. \eqref{eq8}, with the surjective arrow being the identity map on $\mathbbm{1}$. 

\subsection{} Fix an integer $p$. After applying a fiber functor $\omega$ to the identification 
\[
\inHom(M/W_pM, W_pM) \ \subset \ W_{-1}\inEnd(M)
\]
we get an identification
\[
Hom(\omega M/\omega W_pM , \omega W_pM) \ \subset \ W_{-1}End(\omega M),
\]
which thinks of $f : \omega M/W_pM\longrightarrow \omega W_pM $ as the composition Eq. \eqref{eq7}. This way, 
\begin{equation}\label{eq107}
Hom(\omega M/\omega W_pM , \omega W_pM)
\end{equation}
becomes an abelian Lie subalgebra of $W_{-1}End(\omega M)$. The exponential map
\[
\exp \ : \ W_{-1}End(\omega M) \ \longrightarrow \ U(M,\omega)(K) \subset GL(\omega M)(K)
\]
is given by the usual exponential series. On the Lie subalgebra Eq. \eqref{eq107}, it is simply given by
\[
\exp(f) \ = \ I+f.
\]

\subsection{} In this subsection we shall introduce certain Lie subalgebras of $\underline{\fu}(M)$ and subgroups of $\mathcal{U}(M,\omega)$ (for any $\omega$) which play a crucial role in the paper. For any integer $p$, let 
\[
\underline{\fu}_p(M) \ := \ \underline{\fu}(M) \cap \inHom(M/W_pM, W_pM)
\]
and for any $\omega$,
\[
\fu_p(M,\omega) \ := \ \omega\, \underline{\fu}_p(M) \ = \ \fu(M,\omega) \cap Hom(\omega M/\omega W_pM, \omega W_pM).
\]
Then $\fu_p(M,\omega)$ is an abelian Lie subalgebra of $\fu(M,\omega)$. 

For any Lie subalgebra $\mathfrak{l}$ of $W_{-1}End(\omega M)$, we denote the subgroup of $U(M,\omega)$ whose Lie algebra is $\mathfrak{l}$ by $\mathcal{exp}(\mathfrak{l})$ (thus $\mathcal{exp}(\mathfrak{l})(K)=\exp(\mathfrak{l})$). Set
\begin{align*}
\mathcal{U}_p(M,\omega) \ :=& \ \mathcal{exp}(\fu_p(M,\omega)) \\
\ =& \ \mathcal{U}(M,\omega) \ \cap \ \mathcal{exp}(Hom(\omega M/\omega W_pM, \omega W_pM))\\
\ =& \ \mathcal{G}(M,\omega) \ \cap \ \mathcal{exp}(Hom(\omega M/\omega W_pM, \omega W_pM)).
\end{align*}
This is an abelian unipotent subgroup of $\mathcal{U}(M,\omega)$.

\begin{lemma}\label{lem U_p}
$\mathcal{U}_p(M,\omega)$ is the kernel of the restriction homomorphism
\[
\mathcal{G}(M,\omega) \ \longrightarrow \ \mathcal{G}(W_pM\oplus (M/W_pM),\omega)
\]
(induced by $\langle W_pM\oplus (M/W_pM)\rangle^{\otimes} \subset \langle M\rangle^{\otimes}$).
\end{lemma}

\begin{proof}
Tentatively, let us refer to the kernel of the homomorphism given in the statement of the lemma as $U'$. It is clear that $U'$ is contained in $\mathcal{U}(M,\omega)$. In particular, $U'$ is also unipotent and thus it is enough to show that $U'$ and $\mathcal{U}_p(M,\omega)$ have the same $K$-valued points. We have
\[
\mathcal{U}_p(M,\omega)(K) \ = \ \mathcal{G}(M,\omega)(K) \ \cap \ \exp(Hom(\omega M/\omega W_pM, \omega W_pM)).
\]
Let $\sigma\in \mathcal{G}(M,\omega)(K)$. Then $\sigma\in U'(K)$ if and only if $\sigma_{W_pM}=I$ and $\sigma_{M/W_pM}=I$. Under the identification $\mathcal{G}(M,\omega) \subset G(M,\omega)$ (via $\sigma\mapsto \sigma_M$), $\sigma_{W_pM}$ is simply the restriction $\sigma |_{\omega W_pM}$ of $\sigma$ to $\omega W_pM$, and $\sigma_{M/W_pM}$ is the map $\overline{\sigma}$ that $\sigma$, as an element of the parabolic subgroup $G(M,\omega)$, induces on $\omega M/\omega W_pM$ (given by $\overline{\sigma}(v+\omega W_p)= \sigma(v)+\omega W_p$, where $v\in \omega M$). On recalling that
\[
\exp(Hom(\omega M/\omega W_pM, \omega W_pM)) \ = \ I+Hom(\omega M/\omega W_pM, \omega W_pM), 
\]
it is easy to see that the subgroup of $G(M,\omega)(K)$ which acts as identity on both $\omega W_pM$ and $\omega M/\omega W_pM$ is
\[
\exp(Hom(\omega M/\omega W_pM, \omega W_pM)).
\]
The claim follows.
\end{proof}

\begin{rem}
Our examples in Section \ref{counterexample to Thm3 without IA hypothesis} (also see Remark (3) therein) show that in general, $\sum\limits_p \underline{\fu}_p(M)$ may be strictly smaller than $\underline{\fu}(M)$. It is however true that if every $\mathcal{E}_p/\underline{\fu}$ splits, then $\underline{\fu}(M)=\sum\limits_p \underline{\fu}_p(M)$. See Remark (2) at the end of Section \ref{sec: 5.1}.
\end{rem}

\subsection{}\label{on statement star} Let us recall Statement ($\ast$) from the Introduction:

($\ast$) : {\it For any subobject $A$ of $\inHom(M/W_pM,W_pM)$, we have $\underline{\fu}_p(M)\subset A$ if and only if the quotient
\[
\inHom(M/W_pM, M)^\dagger \, / A
\]
belongs to the subcategory Eq. \eqref{eq106}.}

This follows from Thoerem 3.3.1 of \cite{EM21}, which itself is obtained by a slight modification of Hardioun's argument for her Theorem 2 in \cite{Har11} (see also \cite{Har06}). In the interest of keeping the paper more self-contained, let us briefly recall the argument here: The statement is trivial if $M/W_pM=0$ so we may assume otherwise. Let $A$ be a subobject of $\inHom(M/W_pM,W_pM)$ and $\omega$ a fiber functor. In view of Section \ref{section criterion for objects being in subcategories in terms of the action of the kernel of the restriction map} and Lemma \ref{lem U_p}, the quotient 
\[
\inHom(M/W_pM, M)^\dagger \, / A
\]
belongs to the subcategory Eq. \eqref{eq106} if and only if $\mathcal{U}_p(M,\omega)$ acts trivially on 
\begin{equation}\label{eq102}
\omega(\inHom(M/W_pM, M)^\dagger \, / A) \ = \ \omega\, \inHom(M/W_pM, M)^\dagger \, / \omega\, A.
\end{equation}
Choose a section of the natural surjection $\omega M\longrightarrow \omega M/\omega W_pM$ to identify
\[
\omega M \ = \ \omega W_pM \ \oplus \ \omega M/\omega W_pM
\]
(as vector spaces). This also gives a decomposition of $\omega \, \inHom(M/W_pM, M)$. In view of the sequence Eq. \eqref{eq8} and on noting that $\inHom(M/W_pM, W_pM)$ belongs to Eq. \eqref{eq106}, the group $\mathcal{U}_p(M,\omega)$ acts trivially on Eq. \eqref{eq102} if and only if it (or equivalently, $\mathcal{U}_p(M,\omega)(K)$) fixes the image of the element 
\begin{align*}
(0, I) \ \in \ Hom(\omega M/\omega W_pM, \omega M)^\dagger \ \subset & \ Hom(\omega M/\omega W_pM, \omega M) \\ = & \ Hom(\omega M/\omega W_pM, \omega W_pM) \\
&  \oplus End(\omega M/\omega W_pM)
\end{align*}
in Eq. \eqref{eq102}. Now one easily sees that the latter statement is equivalent to $\omega A$ containing $\fu_p(M,\omega)$. (See \cite[Section 3.3]{EM21} for details; note that $L$, $N$, and $\mathcal{U}(M)$ of {\it loc. cit.} would be respectively $W_pM$, $M/W_pM$, and $\mathcal{U}_p(M,\omega)$ here.)

\subsection{} Statement ($\ast$) can be reformulated in the language of extensions originating from subcategories of $\mathbf{T}$ (see Section \ref{section definition of extensions coming from subcategories and lemmas}) as follows:

\begin{thm}\label{thm 1} Let $A$ be a subobject of $\inHom(M/W_pM, W_pM)$. Then the extension $\mathcal{E}_p(M)/A$, viewed as an extension of $\mathbbm{1}$ by $W_{-1}\inEnd(M)/A$ or $\inHom(M/W_pM, W_pM)/A$, originates from the subcategory Eq. \eqref{eq106} if and only if $A$ contains $\underline{\fu}_p(M)$. 
\end{thm}

In other words, $\underline{\fu}_p(M)$ is the smallest subobject of $\inHom(M/W_pM, W_pM)$ such that the extension $\mathcal{E}_p(M)/\underline{\fu}_p(M)$, viewed as an extension of $\mathbbm{1}$ by $W_{-1}\inEnd(M)/ \underline{\fu}_p(M)$ or by $\inHom(M/W_pM, W_pM)/ \underline{\fu}_p(M)$, originates from the subcategory Eq. \eqref{eq106}.

\begin{proof}[Proof of Theorem \ref{thm 1}]
Let $A$ be a subobject of $\inHom(M/W_pM, W_pM)$. By Lemma \ref{lem 3 extensions originating from subcategories} (also see the remark after the same lemma), and in view of the facts (1) that the extension $\mathcal{E}_p(M)/A$ of $\mathbbm{1}$ by $W_{-1}\inEnd(M)/A$ is the image of its namesake as an extension of $\mathbbm{1}$ by $\inHom(M/W_pM, W_pM)/A$ under the obvious pushforward map 
\[
Ext^1(\mathbbm{1}, \inHom(M/W_pM, W_pM)/A) \ \longrightarrow \ Ext^1(\mathbbm{1}, W_{-1}\inEnd(M)/A),
\]
and (2) that (by weight considerations) there are no nonzero morphisms from $\mathbbm{1}$ to objects of weight $<0$, the following statements are equivalent for any full Tannakian subcategorty $\mathbf{S}$ of $\mathbf{T}$ which is closed under taking subobjects:
\begin{itemize}
\item[(i)] The extension $\mathcal{E}_p(M)/A$, viewed as an element of 
\[Ext^1(\mathbbm{1}, W_{-1}\inEnd(M)/A),\] 
originates from $\mathbf{S}$.
\item[(ii)] The extension $\mathcal{E}_p(M)/A $, viewed as an element of
\[Ext^1(\mathbbm{1}, \inHom(M/W_pM, W_pM)/A),\] 
originates from $\mathbf{S}$.
\end{itemize}
In view of Lemma \ref{lem 2 extensions originating from subcategories} and on recalling the explicit description of 
\[\mathcal{E}_p(M) \ \in \ Ext^1(\mathbbm{1}, \inHom(M/W_pM, W_pM))\]
from Section \ref{explicit description of E_p}, Statement (ii) with $\mathbf{S}$ taken to be the subcategory Eq. \eqref{eq106} is equivalent to the following statement:
\begin{itemize}
\item[(iii)] The object 
\[\inHom(M/W_pM, M)^\dagger /A\] 
belongs to Eq. \eqref{eq106}.
\end{itemize}
Thus Theorem \ref{thm 1} is equivalent to Statement ($\ast$) of Section \ref{on statement star} (or the Introduction).
\end{proof}

\section{Extension classes and subgroups of the fundamental group - Part II}\label{section: part 1.2}

In the previous section we saw that $\underline{\fu}_p(M)$ is the smallest subobject of $\inHom(M/W_pM, W_pM)$ such that the extension $\mathcal{E}_p(M)/\underline{\fu}_p(M)$ originates from the subcategory Eq. \eqref{eq106}. Our goal in this section is to give criteria under which the subcategory Eq. \eqref{eq106} in this statement can be replaced by smaller subcategories. Of particular interest will be when we can replace Eq. \eqref{eq106} with a semisimple category.

\subsection{}\label{sec: 5.1} Let us first make an observation regarding the pushforwards of the extension $\mathcal{E}_p(M)$. Recall that we are using the same notation for 
\[
\mathcal{E}_p(M) \ \in Ext^1(\mathbbm{1}, \inHom(M/W_pM, W_pM))
\]
and its image in $Ext^1(\mathbbm{1}, W_{-1}\inEnd(M))$ under the pushforward map Eq. \eqref{eq3}. 

\begin{lemma}\label{lem equivalent statements to E_p/u_p originaring from S}
Let $\mathbf{S}$ be a full Tannakian subcategory of $\mathbf{T}$ closed under taking subobjects. Then the following statements are equivalent:
\begin{itemize}
\item[(i)] The extension 
\[\mathcal{E}_p(M)/\underline{\fu}_p(M) \ \in \  Ext^1(\mathbbm{1}, \inHom(M/W_pM, W_pM)/\underline{\fu}_p(M))\] 
originates from $\mathbf{S}$.
\item[(ii)] The extension 
\[\mathcal{E}_p(M)/\underline{\fu}_p(M) \ \in \  Ext^1(\mathbbm{1}, W_{-1}\inEnd(M)/\underline{\fu}_p(M))\] 
originates from $\mathbf{S}$.
\item[(iii)] The extension 
\[\mathcal{E}_p(M)/\underline{\fu}(M) \ \in \  Ext^1(\mathbbm{1}, W_{-1}\inEnd(M)/\underline{\fu}(M))\] 
originates from $\mathbf{S}$.
\end{itemize}
\end{lemma}

\begin{proof}
That (i) implies (ii) and (ii) implies (iii) is clear, as under the obvious maps the extension in (i) pushes forward to the extension in (ii) and then to the one in (iii) (in fact, we already observed the equivalence of (i) and (ii) in the proof of Theorem \ref{thm 1}). That (iii) implies (i) follows similarly as in the proof of Theorem \ref{thm 1} from Lemma \ref{lem 3 extensions originating from subcategories} on recalling that
\[
\underline{\fu}(M) \, \cap \, \inHom(M/W_pM, W_pM) \ = \ \underline{\fu}_p(M)
\]
(so that the obvious map  
\[
\inHom(M/W_pM, W_pM)/\underline{\fu}_p(M) \ \longrightarrow \ W_{-1}\inEnd(M)/\underline{\fu}(M)
\]
is injective).
\end{proof}

\begin{rem}
\begin{itemize}
\item[(1)] In particular, by taking $\mathbf{S}$ to be the semisimple subcategory $\langle \mathbbm{1}\rangle^{\otimes}$ we see that the three extensions in the lemma split at the same time.
\item[(2)] The lemma together with Deligne's Theorem \ref{Deligne thm} implies that if every $\mathcal{E}_p(M)/\underline{\fu}(M)$ splits (i.e. if every $\mathcal{E}_p(M)$ is in the image of Eq. \eqref{eq101}), then $\underline{\fu}(M)=\sum\limits_p  \underline{\fu}_p(M)$. Indeed, let us tentatively set $\underline{\fu}'=\sum\limits_p  \underline{\fu}_p(M)$. If $\mathcal{E}_p(M)/\underline{\fu}(M)$ splits for every $p$, then so does $\mathcal{E}_p(M)/\underline{\fu}_p(M)$ and hence $\mathcal{E}_p(M)/\underline{\fu}'$ (the latter as an extension of $\mathbbm{1}$ by $W_{-1}\inEnd(M)/\underline{\fu}'$). It follows that $\mathcal{E}(M)/\underline{\fu}'$ splits, so that by Deligne's theorem $\underline{\fu}(M)\subset \underline{\fu}'$.
\end{itemize}
\end{rem}

\subsection{}\label{par: IA axioms} For any integers $p$ and $q$ with $q\leq p$, define
\begin{align*}
J_1^{\{p,q\}} \ :=& \ \{(i,j)\in \ZZ^2 :i\leq p <j\}\\
J_2^{\{p,q\}} \ :=& \ \{(i,j)\in \ZZ^2 : i<j \ \text{and} \ (q <j \leq p \ \text{or} \ i> p )\}.
\end{align*}

Figure 1 shows the two sets. In the figure, the axes are oriented according to the standard labelling of entries of a matrix (the pair $(i,j)$ is placed where the entry $ij$ of a matrix sits). 

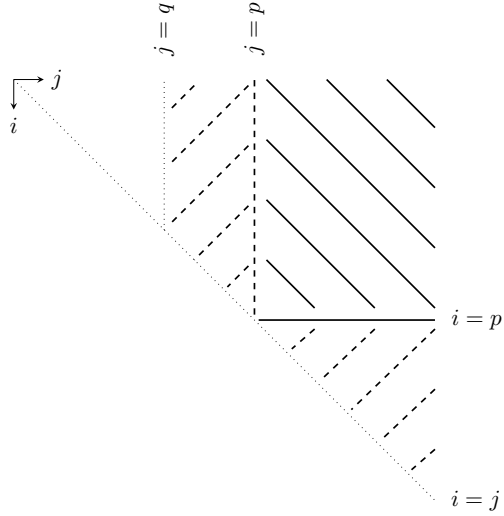
\begin{figure}\label{fig1}
\scalebox{.8}{
\begin{tikzpicture}
\draw[->] (0,7) -- (0,6.5) node[anchor=north] {$i$};
\draw[->] (0,7) -- (0.5,7) node[anchor=west] {$j$};
\draw[dotted] (0,7) -- (7,0) node[anchor=west]{ \ $i=j$};
\draw[dotted] (2.5,4.5) -- (2.5,7) node[anchor=south]{\rotatebox{90}{ \ $j=q$}};
\draw[thick, dashed] (4,3.09) -- (4,7) node[anchor=south]{\rotatebox{90}{ \ $j=p$}};
\draw[thick] (4.07,3) -- (7,3) node[anchor=west]{ \ $i=p$}; 

\draw[thick] (5,3.2) -- (4.2,4);
\draw[thick] (6,3.2) -- (4.2,5);
\draw[thick] (7,3.2) -- (4.2,6);
\draw[thick] (7,4.2) -- (4.2,7);
\draw[thick] (7,5.2) -- (5.2,7);
\draw[thick] (7,6.2) -- (6.2,7);

\draw[thick, dashed] (3.9,3.9) -- (3.55,3.55);
\draw[thick, dashed] (3.9,4.9) -- (3,4);
\draw[thick, dashed] (3.9,5.9) -- (2.6,4.6);
\draw[thick, dashed] (3.9,6.9) -- (2.6,5.6);
\draw[thick, dashed] (3,6.9) -- (2.6,6.5);

\draw[thick, dashed] (5, 2.85) -- (4.62,2.52);
\draw[thick, dashed] (6, 2.85) -- (5.1,2);
\draw[thick, dashed] (7, 2.85) -- (5.61,1.51);
\draw[thick, dashed] (7, 1.85) -- (6.11,1.01);
\draw[thick, dashed] (7, 0.85) -- (6.61,0.51);
\end{tikzpicture}}
\caption{The set of lattice points in the region marked by solid (resp. thick dashed) lines is $J_1^{\{p,q\}}$ (resp. $J_2^{\{p,q\}}$).}
\end{figure}

We consider the following {\it independence axioms} for an object $M$ of $\mathbf{T}$:

\begin{itemize}
\item $(IA1)_{\{p,q\}}$: The two objects 
\[
\bigoplus_{(i,j)\in J_1^{\{p,q\}}} \inHom(Gr^W_{j}M,Gr^W_{i}M)  \ \ \ \  \ \text{and} \ \ \ \ \bigoplus_{(i,j)\in J_2^{\{p,q\}}} \inHom(Gr^W_{j}M,Gr^W_{i}M)
\]
have no nonzero isomorphic subobjects. (Note that if $q'\leq q\leq p$, then $(IA1)_{\{p,q'\}}$ implies $(IA1)_{\{p,q\}}$.

\item $(IA2)_{\{p,q\}}$: The two sets
\[
J_1^{\{p,q\}}(M) \ := \ \{i-j : (i,j)\in J_1^{\{p,q\}}, \, Gr^W_iM\neq 0, \,  Gr^W_jM\neq 0\}
\]
and 
\[
J_2^{\{p,q\}}(M) \ := \ \{i-j : (i,j)\in J_2^{\{p,q\}}, \, Gr^W_iM\neq 0, \,  Gr^W_jM\neq 0\}
\]
are disjoint. (Note that $J_1^{\{p,q\}}(M)$ and $J_2^{\{p,q\}}(M)$ are respectively the set of weights of the two object in $(IA1)_{\{p,q\}}$ above.)
\item $(IA3)$: The numbers 
\[
i-j \hspace{.5in}(i<j, \ Gr^W_iM\neq 0, \,  Gr^W_jM\neq 0)
\]
are all distinct. (Equivalently, if $M$ has $n$ distinct weights, then $W_{-1}\inEnd(M)$ has ${n \choose 2}$ distinct weights.)
\end{itemize}

It is clear that $(IA2)_{\{p,q\}}$ implies $(IA1)_{\{p,q\}}$, and $(IA3)$ implies $(IA2)_{\{p,q\}}$ for every $p$ and $q$. Also note that whether or not $M$ satisfies any of these axioms only depends on $Gr^WM$.

\subsection{} We can now state the main result of this part of the paper:

\begin{thm}\label{thm 2}
Let $q\leq p$. Consider the following statements:
\begin{itemize}
\item[(i)] $M$ satisfies $(IA1)_{\{p,q\}}$ and $Gr^WM$ is semisimple ( = completely reducible).
\item[(ii)] $M$ satisfies $(IA2)_{\{p,q\}}$. 
\end{itemize}
If either statement holds, then the extension $\mathcal{E}_p(M)/\underline{\fu}_p(M)$ originates from the subcategory $\langle W_qM , Gr^WM\rangle^{\otimes}$.
\end{thm}

The proof of Theorem \ref{thm 2} shall be given in the Sections \ref{introducing U_>=q} - \ref{proof of thm 2} below. Here we consider some consequences of the theorem:

\begin{enumerate}
\item Since $q\leq p$, the subcategegory $\langle W_qM , Gr^WM\rangle^{\otimes}$ is contained in the subcategory $\langle W_pM , M/W_pM\rangle^{\otimes}$. Thus combining Thoerems \ref{thm 1} and \ref{thm 2} we get the following refinement of Theorem \ref{thm 1}: if Statements (i) or (ii) above hold for some $q\leq p$, then $\underline{\fu}_p(M)$ is the smallest subobject of $\inHom(M/W_pM, W_pM)$ such that $\mathcal{E}_p(M)/\underline{\fu}_p(M)$ originates from $\langle W_qM , Gr^WM\rangle^{\otimes}$.

\item Perhaps the most interesting application of Theorem \ref{thm 2} is in the following scenario: Fix $p$. Suppose $Gr^WM$ is semisimple; for instance, this will be the case if $\mathbf{T}$ is a category of motives, or if $\mathbf{T}$ is the category of mixed Hodge structures and $Gr^WM$ is polarizable. Suppose $M$ satisfies $(IA1)_{\{p,q\}}$ for all $q\leq p$ (this holds for instance, if $M$ satisfies $(IA3)$). Then $\underline{\fu}_p(M)$ is the smallest subobject of $\inHom(M/W_pM, W_pM)$ such that $\mathcal{E}_p(M)/\underline{\fu}_p(M)$ originates from the semisimple subcategory $\langle Gr^WM\rangle^{\otimes}$, i.e. splits. In particular, $\mathcal{E}_p(M)/\underline{\fu}(M)$ splits. For future referencing, we record this as a corollary:
\end{enumerate}

\begin{cor}\label{cor 1 of thm 2} Fix $p$. Suppose $Gr^WM$ is semisimple and that $M$ satisfies $(IA1)_{\{p,q\}}$ for all $q\leq p$. Then $\underline{\fu}_p(M)$ is the smallest subobject of $\inHom(M/W_pM, W_pM)$ such that $\mathcal{E}_p(M)/\underline{\fu}_p(M)$ splits. In particular, 
\[
\mathcal{E}_p(M)/\underline{\fu}(M)
\]
splits.
\end{cor}

As a special case, we obtain:
\begin{cor}\label{cor 2 of thm 2} If $Gr^WM$ is semisimple and $(IA3)$ holds, then for every $p$ the extension $\mathcal{E}_p(M)/\underline{\fu}(M)$ splits.
\end{cor}

\begin{rem}
Recall that by Deligne's Theorem \ref{Deligne thm}, the extension 
\[\sum\limits_p \mathcal{E}_p(M)/\underline{\fu}(M)\] 
splits. As we pointed out earlier, in general, the individual extensions $\mathcal{E}_p(M)/\underline{\fu}(M)$ may not split (see Section \ref{counterexample to Thm3 without IA hypothesis} and Remark (2) therein for examples). The above results give sufficient conditions for when an individual $\mathcal{E}_p(M)/\underline{\fu}(M)$ splits.
\end{rem}

\subsection{}\label{introducing U_>=q} From this point until the end of Section \ref{proof of thm 2} our goal is to prove Theorem \ref{thm 2}. Given any fiber functor $\omega$, let $\mathcal{U}_{\geq q}(M,\omega)$ be the kernel of the surjection 
\[
\mathcal{G}(M,\omega) \ \longrightarrow \ \mathcal{G}(W_qM\oplus Gr^WM,\omega)
\]
induced by the inclusion $\langle W_qM\oplus Gr^WM\rangle^{\otimes}\subset \langle M\rangle^{\otimes}$. Then $\mathcal{U}_{\geq q}(M,\omega)$ is the subgroup of  $\mathcal{U}(M,\omega)$ which acts trivially on $\omega W_qM$. Let $U_{\geq q}(M,\omega)$ be the subgroup of $GL(\omega M)$ consisting of the elements which fix the weight filtration, and act trivially on $Gr^W\omega M$ and $\omega W_qM$:
\[
U_{\geq q}(M,\omega) \ := \ \{\sigma\in U(M,\omega): \sigma |_{\omega W_qM}=I\}. 
\]
Then
\[
\mathcal{U}_{\geq q}(M,\omega) \ = \ \mathcal{U}(M,\omega) \cap U_{\geq q}(M,\omega).
\]
We have
\[
Lie (U_{\geq q}(M,\omega)) \ = \ Hom(\omega M/\omega W_qM, \omega M)\cap W_{-1}End(\omega M),
\]
where $Hom(\omega M/\omega W_qM, \omega M)$ is identified as the subspace of $End(\omega M)$ consisting of the elements which vanish on $\omega W_qM$. Then
\[
\fu_{\geq q}(M,\omega) \ := \ Lie (\mathcal{U}_{\geq q}(M,\omega)) \ = \ \fu(M,\omega) \, \cap \, Hom(\omega M/\omega W_qM, \omega M).
\]
Finally, set 
\[
\underline{\fu}_{\geq q}(M) \ := \  \underline{\fu}(M) \, \cap \, \inHom(M/W_qM, M),
\]
where the intersection means fibered product over $\inEnd(M)$. Here 
\[\inHom(M/W_qM, M)\] 
is thought of as a subobject of $\inEnd(M)$ via the obvious injection induced by the quotient map $M\longrightarrow M/W_qM$ (note that this is compatible with the previous identification of $Hom(\omega M/\omega W_qM, \omega M)$ as a subspace of $End(\omega M)$). We then have
\[
\fu_{\geq q}(M,\omega) \ = \ \omega\, \underline{\fu}_{\geq q}(M).
\]

\subsection{}\label{how IA axioms are relevant} Identifying
\begin{equation}\label{eq13}
Gr^W \inEnd(M) \ = \ \inEnd(Gr^WM) \ = \ \bigoplus_{i,j} \, \inHom(Gr^W_jM,Gr^W_iM),
\end{equation}
we have
\[Gr^W W_{-1}\inEnd(M) \ = \ \bigoplus_{\stackrel{i,j}{i<j}} \, \inHom(Gr^W_jM,Gr^W_iM).\]
Then for every $q$,
\begin{align}
Gr^W\underline{\fu}_{\geq q}(M) \, &\subset \, Gr^W \inHom(M/W_qM, M) \, \cap \, Gr^W W_{-1}\inEnd(M) \notag \\
&= \ \bigoplus_{\stackrel{i,j}{i,q< j}} \inHom(Gr^W_jM,Gr^W_iM).\label{eq15}
\end{align}

The following lemma is the only place in the proof of Theorem \ref{thm 2} that conditions (i) and (ii) of the theorem play a part.

\begin{lemma}\label{decomposition of Gr u_>=p}
Let $q\leq p$. Suppose Statement (i) or (ii) of Theorem \ref{thm 2} holds. Then $Gr^W\underline{\fu}_{\geq q}(M)$ decomposes as the direct sum of 
\[Gr^W\underline{\fu}_{\geq q}(M) \ \  \cap \ \bigoplus_{(i,j)\in J_1^{\{p,q\}}} \inHom(Gr^W_{j}M,Gr^W_{i}M)\]
and 
\[Gr^W\underline{\fu}_{\geq q}(M) \  \ \cap \ \bigoplus_{(i,j)\in J_2^{\{p,q\}}} \inHom(Gr^W_{j}M,Gr^W_{i}M).\]
\end{lemma}

\begin{proof}
The direct sum in Eq. \eqref{eq15} is over all pairs $(i,j)$ in $J_1^{\{p,q\}}\sqcup J_2^{\{p,q\}}$, so that we can rewrite Eq. \eqref{eq15} as
\begin{align*}
Gr^W\underline{\fu}_{\geq q}(M) \ \subset& \ \overbrace{\bigoplus_{(i,j)\in J_1^{\{p,q\}}} \inHom(Gr^W_{j}M,Gr^W_{i}M)}^{\text{(I)}} \\ 
 & \ \ \oplus \ \overbrace{\bigoplus_{(i,j)\in J_2^{\{p,q\}}} \inHom(Gr^W_{j}M,Gr^W_{i}M)}^{\text{(II)}}.
\end{align*}
First suppose $Gr^WM$ is semisimple and $M$ satisfies $(IA1)_{\{p,q\}}$. Then the object $Gr^W\underline{\fu}_{\geq q}(M)$ (living in the semisimple category $\langle Gr^WM\rangle^{\otimes}$) is a direct sum of simple objects. By $(IA1)_{\{p,q\}}$, each simple direct factor either lives in (I) or (II).

On the other hand, if $(IA2)_{\{p,q\}}$ holds, then each nonzero graded component $Gr_n^W\underline{\fu}_{\geq q}(M)$ must live in (I) or (II) (whichever has a nonzero weight $n$ part). 
\end{proof}

\subsection{}\label{proof of thm 2} We are ready to give the proof of Theorem \ref{thm 2}. We may assume that $M/W_pM$ is not zero. Consider $\mathcal{E}_p(M)$ as an extension of the unit object by $\inHom(M/W_pM, W_pM)$, given by Eq. \eqref{eq8}.  In view of Section \ref{introducing U_>=q} and Lemma \ref{lem 1 extensions originating from subcategories}, it is enough to check right exactness of the sequence obtained by applying $\mathcal{U}_{\geq q}(M,\omega)$-invariance to $\omega(\mathcal{E}_p(M)/\underline{\fu}_p(M))$ for a suitably chosen fiber functor $\omega$. Let $\omega_0$ be an arbitrary fiber functor. We shall take the composition
\[
\omega^{gr} \  : \ \mathbf{T} \ \stackrel{Gr^W}{\longrightarrow} \ \mathbf{T} \ \stackrel{\omega_0}{\longrightarrow} \   \mathbf{Mod}_K.
\]
as our fiber functor.

Via the identification
\[
\inHom \, (M/W_pM, M) \ \subset \ \inEnd \, (M), 
\]
we think of the image under $\omega^{gr}$ of every subobject of $\inHom(M/W_pM, M)$ as a subspace of $\omega^{gr} \inEnd(M)$. Throughout, we shall write the elements of 
\begin{align*}
\omega^{gr} \inEnd \, (M) \ = \ End \, (\omega^{gr}M) \ =& \ End \, (\bigoplus_n \omega_0 Gr^W_nM) \\
 =& \ \bigoplus\limits_{i,j}Hom \, (\omega_0 Gr^W_jM, \omega_0 Gr^W_iM)
\end{align*}
as 2 by 2 block matrices with rows (resp. columns) broken up as $\{i \, : \, i\leq p\}\, \cup \, \{i \, : \, i>p\}$ (resp. the same with $j$ replacing $i$). Then an element
\[
f \ \in \ \omega^{gr} \inHom(M/W_pM, M)^\dagger \ = \ Hom(\omega^{gr}(M/W_pM), \omega^{gr}M)^\dagger
\]
looks like
\[
\begin{pmatrix}
0 & \ast\\
0 & \lambda(f) I
\end{pmatrix}.
\]
The surjective arrow 
\[
Hom(\omega^{gr}(M/W_pM), \omega^{gr}M)^\dagger \ \longrightarrow \ K
\]
in $\omega^{gr}\mathcal{E}_p$ sends $f$ to $\lambda(f)$. 

Consider the element
\[
f_0 \ = \ \begin{pmatrix}
0 & 0\\
0 &  I
\end{pmatrix} \ \in \ Hom(\omega^{gr}(M/W_pM), \omega^{gr}M)^\dagger .
\]
We will show that if Conditions (i) or (ii) of Theorem \ref{thm 2} hold (and $q\leq p$), then the element $f_0+\omega^{gr}\underline{\fu}_p(M)$ of 
\[
\frac{Hom(\omega^{gr}(M/W_pM), \omega^{gr}M)^\dagger}{\omega^{gr}\underline{\fu}_p(M)}
\]
is fixed by $\mathcal{U}_{\geq q}(M,\omega^{gr})$; this proves surjectivity of 
\[
\left(\frac{Hom(\omega^{gr}(M/W_pM), \omega^{gr}M)^\dagger}{\omega^{gr}\underline{\fu}_p(M)}\right)^{\mathcal{U}_{\geq q}(M,\omega^{gr})} \ \longrightarrow \ K
\]
and hence the theorem. Since $\mathcal{U}_{\geq q}(M,\omega^{gr})$ is unipotent, it is enough to verify that $f_0+\omega^{gr}\underline{\fu}_p(M)$ is fixed by every $\sigma\in \mathcal{U}_{\geq q}(M,\omega^{gr})(K)\subset GL(\omega^{gr} M)$. Given such $\sigma$, we must show that
\begin{equation}\label{eq14}
\sigma f_0 \sigma^{-1} \, - \, f_0 \ \in \ \omega^{gr}\underline{\fu}_p(M) \, = \, \fu_p(M,\omega^{gr}). 
\end{equation}

Writing
\[
\sigma \ = \ \begin{pmatrix}
\sigma_1 & A\\
0 & \sigma_2
\end{pmatrix},
\]
we have 
\[
\log(\sigma) \ = \ \begin{pmatrix}
\log \sigma_1 & \ast \\
0 & \log \sigma_2
\end{pmatrix} \ \in \ \fu_{\geq q}(M,\omega^{gr}) \ = \ \omega_0 Gr^W\underline{\fu}_{\geq q}(M).
\]
Applying $\omega_0$ to the decomposition of $Gr^W\underline{\fu}_{\geq q}(M)$ given in Lemma \ref{decomposition of Gr u_>=p} it follows that 
\[
\begin{pmatrix}
\log \sigma_1 &  0 \\
0 & \log \sigma_2
\end{pmatrix} \ \in \ \fu_{\geq q}(M,\omega^{gr}),
\]
so that
\[
\delta  \ := \ \begin{pmatrix}
\sigma_1 &  0 \\
0 & \sigma_2
\end{pmatrix} \ \in \ \mathcal{U}_{\geq q}(M,\omega^{gr})(K).
\]
We thus have
\[
\sigma f_0 \sigma^{-1} \, - \, f_0 \ = \ \begin{pmatrix}
0& A\sigma^{-1}_2\\
0 & 0
\end{pmatrix} \ = \ \log(\sigma\delta^{-1} ) \ \in \ \fu_{\geq q}(M,\omega^{gr}),
\]
as desired. (Note that if an element of the form $\begin{pmatrix} 0 & \ast \\ 0 & 0\end{pmatrix}$ is in $\fu(M,\omega^{gr})$, then it will actually be in $\fu_{p}(M,\omega^{gr})$.)

\subsection{} We end this section with a variant of Theorem \ref{thm 2} for $q>p$, which again gives a sufficient condition to guarantee that $\mathcal{E}_p(M)/\underline{\fu}_p(M)$ originates from the category $\langle W_qM , Gr^WM\rangle^{\otimes}$. Note that when $q>p$, then none of the two categories $\langle W_qM , Gr^WM\rangle^{\otimes}$ and $\langle W_pM , M/W_pM\rangle^{\otimes}$ necessarily contains the other\footnote{The content of this subsection will not be used anywhere else in the paper.}.
 
For $q>p$, consider the following three sets:
\begin{align*}
{J'}_1^{\{p,q\}} \ :=& \ \{(i,j)\in \ZZ^2 : \, i\leq p, \ j>q \}\\
{J'}_2^{\{p,q\}} \ :=& \ \{(i,j)\in \ZZ^2 : \,  p< i \leq q <j \}\\
{J'}_3^{\{p,q\}} \ :=& \ \{(i,j)\in \ZZ^2 : \, q< i<j \}.
\end{align*}
Say an object $M$ of $\mathbf{T}$ satisfies $(IA1')_{\{p,q\}}$ if the objects
\[
\bigoplus_{(i,j)\in {J'}_k^{\{p,q\}}} \inHom(Gr^W_jM, Gr^W_iM)
\]
for $k=1,2,3$ have no nonzero isomorphic subobjects. We say $M$ satisfies $(IA2')_{\{p,q\}}$ if the sets of weights of these objects are disjoint. Then $(IA2')_{\{p,q\}}$ implies $(IA1')_{\{p,q\}}$, and $(IA3)$ implies $(IA2')_{\{p,q\}}$ for every $p,q$. 

\begin{thm}\label{thm2'}
Let $q>p$. Suppose one of the following statements holds:
\begin{itemize}
\item[(i)] $Gr^WM$ is semisimple and $M$ satisfies $(IA1')_{\{p,q\}}$.
\item[(ii)] $M$ satisfies $(IA2')_{\{p,q\}}$.
\end{itemize}
Then the extension $\mathcal{E}_p(M)/\underline{\fu}_p(M)$ originates from $\langle W_qM , Gr^WM\rangle^{\otimes}$. 
\end{thm}

\begin{proof}
The proof is similar to the proof of Theorem \ref{thm 2}. Note that the pairs $(i,j)$ appearing in Eq. \eqref{eq15} are those in ${J'}_1^{\{p,q\}}\cup {J'}_2^{\{p,q\}}\cup {J'}_3^{\{p,q\}}$. Similar to Lemma \ref{decomposition of Gr u_>=p}, Hypothesis (i) or (ii) above imply that $Gr^W\underline{\fu}_{\geq q}(M)$ is the direct sum of its intersections with the three objects
\begin{equation}\label{eq16}
\bigoplus_{(i,j)\in {J'}_k^{\{p,q\}}} \inHom(Gr^W_jM, Gr^W_iM)
\end{equation}
for $k=1,2,3$. Taking $\omega^{gr}$ and $f_0$ as in the proof of Theorem \ref{thm 2}, we shall show that
for every $\sigma\in \mathcal{U}_{\geq q}(M,\omega^{gr})(K)$, 
\[
\sigma f_0 \sigma^{-1} - f_0 \ \in \ \fu(M,\omega^{gr}).
\]
(it will then automatically be in $\fu_p(M,\omega^{gr})$). Decompose 
\[
\log\sigma \ = \ \tau_1+\tau_2+\tau_3,
\]
where $\tau_k$ is the component in Eq. \eqref{eq16}; each $\tau_k$ is in $\fu_{\geq q}(M,\omega^{gr})$, thanks to Hypothesis (i) or (ii). Writing the elements of $End(\omega^{gr}M)$ as $3\times 3$ block matrices with the rows (resp. columns) broken up as $\{i: i\leq p\} \cup \{i: p<i\leq q\} \cup \{i: i>q\}$ (resp. the same with $j$ replacing $i$), we have
\[
\log(\sigma) \ = \ \begin{pmatrix}
0 & & \tau_1\\
& 0 &\tau_2 \\
 &  & \tau_3
\end{pmatrix}
\]
(with zero missing entries), so that
\[
\sigma \ = \ \begin{pmatrix}
I & & \tau_1(\exp(\tau_3)-1)/\tau_3\\
& I &\tau_2 (\exp(\tau_3)-1)/\tau_3\\
 &  & \exp(\tau_3)
\end{pmatrix}\]
and 
\[\sigma^{-1} \ = \ \begin{pmatrix}
I & & \tau_1(\exp(-\tau_3)-1)/\tau_3\\
& I &\tau_2 (\exp(-\tau_3)-1)/\tau_3\\
 &  & \exp(-\tau_3)
\end{pmatrix},
\]
where for brevity, for a nilpotent map $N$ we have set
\[
(\exp(N)-1)/N \ := \ \sum\limits_{n\geq 0} \, N^n/(n+1)! \, .
\]
Then one calculates
\[
\sigma f_0 \sigma^{-1} - f_0 \ = \ \begin{pmatrix}
0 & & \tau_1(1-\exp(-\tau_3))/\tau_3\\
& 0 & 0\\
 &  & 0
\end{pmatrix}.
\]
This belongs to $\fu(M,\omega^{gr})$ because $\tau_1,\tau_3$ are in the Lie algebra $\fu(M,\omega^{gr})$ and 
\begin{align*}
[\tau_1,\tau_3]  \ &= \ \tau_1\tau_3\\
[[\tau_1,\tau_3], \tau_3] \ &= \ \tau_1\tau_3^2\\
& \ \, \vdots \quad \quad \quad .
\end{align*}
\end{proof}

\section{Motives with large unipotent radicals of motivic Galois groups}\label{sec: part 2}

\subsection{}\label{par: part 2 setup} In this section, unless otherwise indicated, $\mathbf{T}$ is any reasonable Tannakian category of mixed motives in characteristic zero, or the category of mixed Hodge structures. Examples of the former include the (equivalent) categories of mixed motives over a subfield of $\CC$ due to Nori or Ayoub (see \cite{HM17}, \cite{Ay14} and \cite{Ar13}), Voevodsky's category of mixed Tate motives over $\QQ$ (or those over $\ZZ$, etc.), and categories of mixed motives defined via realizations (see Deligne \cite{De89} or Jannsen \cite{Ja90}). See the remark at the end of this subsection for what we exactly need of $\mathbf{T}$. We shall use the word {\it motive} to refer to any object of $\mathbf{T}$ whose weight associated graded is semisimple. Of course, in the case that $\mathbf{T}$ is a reasonable category of mixed motives, this will simply mean an arbitrary object of $\mathbf{T}$. In the case of the category of mixed Hodge structures, this will include (graded-) polarizable objects, and in particular, the Hodge realizations of mixed motives. 

Let $M$ be a motive. We say $\underline{\fu}(M)$ is {\it large} (or that $M$ has a large $\underline{\fu}$) if 
\[
\underline{\fu}(M) \ = \ W_{-1}\inEnd(M).
\]
Similarly, we say $\underline{\fu}_p(M)$ is large if 
\[
\underline{\fu}_p(M) \ = \ \inHom(M/W_pM, W_pM).
\]
Then $\underline{\fu}(M)$ is large if and only if $\underline{\fu}_p(M)$ is large for every $p$. The interest in motives with large $\underline{\fu}$ is partly because of Grothendieck's period conjecture. If $\mathbf{T}$ is a good category of motives over a number field, among the motives with a fixed associated graded, the periods of a motive with large $\underline{\fu}$ should generate a field with the largest possible transcendence degree. We refer the reader to \cite{An04} for a detailed discussion of the Grothendieck's period conjecture.

Our main goal in this section is to use the earlier results of the paper to construct motives with large $\underline{\fu}$ and three weights. We will be particularly interested in motives $M$ with three weights $-2n<p<0$, associated graded isomorphic to 
\[
\QQ(n) \oplus A \oplus \mathbbm{1}
\]
where $A$ is a given pure motive of weight $p$, and such that $\underline{\fu}(M)$ is large. We shall prove a precise classification result for such motives in terms of homological algebra, which completely classifies such motives up to isomorphism when $n\neq -p$ and $Ext^1(\mathbbm{1},\QQ(n))=0$ (e.g. for even $n$ if $\mathbb{T}$ is any reasonable category of motives over $\QQ$). The condition $n\neq -p$ here is an independence axiom (referring to the language of the previous section). See Corollary \ref{cor: motives with three weights and Gr^W isomorphic to Q(n)+A+1, bijective case } for the precise statement of the classification result. As an application, in Section \ref{sec: classification of 3 dimensional mixed Tate motives with large u} we shall consider the case where $A$ is the simple Tate motive $\QQ(k)$ and construct certain interesting mixed Tate motives over $\QQ$.

It turns out that the machinery we shall need works in more generality with little extra effort. So we have decided to develop the results in more generality first and then apply them to the case of motives with three weights. We shall however start with the simplest case below, i.e. motives with only two weights; the observations made in this case will be useful when we deal with more than two weights.

\begin{rem}
Our restriction to the categories of motives and mixed Hodge structures here is for reasons to do with motivation and applications. Unless we explicitly mention otherwise, the discussions can be thought of to take place in the following setting: Take $\mathbf{T}$ to be any Tannakian category equipped with a weight filtration (as in previous sections), and interpret the word ``motive" as an object of $\mathbf{T}$ whose associate graded with respect to the weight filtration is semisimple. In discussions where the Tate objects $\QQ(n)$ make an appearance, we also need to assume that $\mathbf{T}$ contains an object $\QQ(1)$, which is pure of weight -2 and such that the functor $-(1):=-\otimes \QQ(1)$ is invertible with inverse $-\otimes\QQ(-1)$ (and then the $\QQ(n)$ are defined in the usual way from $\QQ(1)$).
\end{rem}

\subsection{}\label{par: case of two weights} We shall use the following terminology: an extension of $\mathbbm{1}$ by an object $L$ is {\it totally nonsplit} if its pushforward to any nonzero quotient of $L$ is nontivial ( = nonsplit); dually, we say an extension of an object $L$ by $\mathbbm{1}$ is totally nonsplit if its pullback to any nonzero subobject of $L$ is nontrivial. Note that if $L$ is simple, then ``totally nonsplit" and ``nonsplit" are equivalent. 

Suppose $M$ is an object with two weights, fitting in a short exact sequence
\begin{equation}\label{eq24}
0 \ \longrightarrow \ L  \ \longrightarrow \ M  \ \longrightarrow \  \mathbbm{1}  \ \longrightarrow \  0,
\end{equation}
where $L$ is a pure motive of weight $p<0$\footnote{Note that this makes $M$ also a motive (as $Gr^WM\simeq L\oplus \mathbbm{1}$ is semisimple).}. Then
\[
W_{-1}\inEnd(M) \ = \ \inHom(\mathbbm{1}, L) \ \cong \ L. 
\]
By Theorem \ref{thm 1} (or Deligne's Theorem \ref{Deligne thm} or Hardouin's \cite[Theorem 2]{Har11}, see also the latter's predecessor, Bertrand's \cite[Theorem 1.1]{Ber01}), $\underline{\fu}(M)$ ( = $\underline{\fu}_p(M)$) is the smallest subobject of $L$ such that the pushforward of the extension Eq. \eqref{eq24} to $Ext^1(\mathbbm{1}, L/\underline{\fu}(M))$ splits. (Indeed, note that via the identification of $\inHom(\mathbbm{1}, L)$ and $L$, the extension $\mathcal{E}_p(M)$ appearing in Theorem \ref{thm 1} is simply Eq. \eqref{eq24}. Also note that the total class $\mathcal{E}(M)$ of $M$ is a nonzero multiple of $\mathcal{E}_p(M)$.) Thus $\underline{\fu}(M)$ is large if and only if Eq. \eqref{eq24} is totally nonsplit. In particular, if $L$ is simple, then
\[
\underline{\fu}(M) \ = \ \begin{cases}
L \hspace{.3in} &\text{if $M$ is not semisimple}\\
0 & \text{if $M$ is semisimple.}\end{cases}
\]

\begin{rem} Let $\mathbf{T}$ be any Tannakian category over $K$ with a weight filtartion. Then for any object $M$ of $\mathbf{T}$ with semisimple  $Gr^WM$ the following statements are equivalent: 
\begin{itemize}
\item[(i)] $\underline{\fu}(M)$ is zero.
\item[(ii)] $M$ is semisimple.
\item[(iii)] $M$ is isomorphic to $Gr^WM$.
\end{itemize}
Indeed, choosing a fiber functor one easily sees $(i)\Rightarrow (ii) \Rightarrow (iii) \Rightarrow (i)$ (note that among these the implication $(i)\Rightarrow (ii)$ is the only one that needs the assumption of semisimplicity of $Gr^WM$). This gives another argument for the characterization of $\underline{\fu}(M)$ given above when $L$ is simple.
\end{rem}

\subsection{}\label{sec: thm 3} In this section we will use the results of Sections \ref{section:  part 1.1} and \ref{section: part 1.2} to give a criterion for a motive to have a large $\underline{\fu}$ in terms of its subobjects and subquotients. 

\begin{thm}\label{thm3}
Let $p<0$ and $M$ be a motive such that
\begin{equation}\label{eq28}
M/W_pM \ \simeq \ \mathbbm{1}, \quad Gr^W_pM \ \neq \ 0
\end{equation}
(so that in particular, $0$ and $p$ are the highest two weights of $M$). Suppose moreover that:
\begin{itemize}
\item[(i)] $\underline{\fu}(W_pM)$ is large, 
\item[(ii)] $\underline{\fu}(M/W_{p-1}M)$ is large, and
\item[(iii)] $M$ satisfies $(IA1)_{\{p,q\}}$ for all $q\leq p$.
\end{itemize}
Then $\underline{\fu}(M)$ is large.
\end{thm}

\begin{proof}
Note that since $M/W_pM$ is pure, for any choice of fiber functor $\omega$, we have $\mathcal{U}_{\geq p}(M,\omega) = \mathcal{U}_{p}(M,\omega)$. (Indeed, if $\sigma$ is in $\mathcal{G}(M,\omega)$, then $\sigma_{Gr^WM}$ and $\sigma_{W_pM}$ are both identity if and only if $\sigma_{Gr^W(M/W_pM)}$ and $\sigma_{W_pM}$ are both identity, and that here $Gr^W(M/W_pM)\simeq M/W_pM$.) Thus the kernel of the natural surjection (or rather, epimorphism)
\[
\underline{\fu}(M)  \ \longrightarrow \ \underline{\fu}(W_pM) 
\]
(induced by the inclusion $\langle W_pM\rangle^{\otimes}\subset \langle M\rangle^{\otimes}$) is $\underline{\fu}_p(M)$. From this (and in light of purity of $M/W_pM$) it follows (e.g. by comparing maximum possible dimensions) that $\underline{\fu}(M)$ is large if and only if $\underline{\fu}(W_pM)$ and $\underline{\fu}_p(M)$ are both large. 

By Corollary \ref{cor 1 of thm 2} (and in view of Hypothesis (iii) and the fact that $M$ is a motive), $\underline{\fu}_p(M)$ is the smallest subobject of $\inHom(M/W_pM, W_pM)$ such that $\mathcal{E}_p(M)/\underline{\fu}_p(M)$ splits. Fix an isomorphism between $M/W_pM$ and $\mathbbm{1}$ to identify the two objects. Then
\[
\underline{\fu}_p(M) \ \subset \ \inHom(M/W_pM, W_pM) \ = \ \inHom(\mathbbm{1}, W_pM) \ \cong \ W_pM.
\]
Via the latter identification, the extension 
\[\mathcal{E}_p(M) \ \in \ Ext^1(\mathbbm{1}, \inHom(\mathbbm{1}, W_pM)) \ = \ Ext^1(\mathbbm{1}, W_pM) \]
is simply the canonical extension
\begin{equation}\label{eq26}
0 \ \longrightarrow \ W_pM \ \longrightarrow \ M  \ \longrightarrow \ \mathbbm{1} \ \longrightarrow \  0
\end{equation}
(where the surjective arrow is the quotient map $M\longrightarrow M/W_pM=\mathbbm{1}$). Let $A$ be any subobject of $W_pM$ such that $\mathcal{E}_p(M)/A$ splits. The goal is to show that $A=W_pM$.

Modding out by $W_{p-1}M$ the extension Eq. \eqref{eq26} pushes forward to
\begin{equation}\label{eq25}
0 \ \longrightarrow \ Gr^W_pM \ \longrightarrow \ M/W_{p-1}M  \ \longrightarrow \ \mathbbm{1} \ \longrightarrow \  0.
\end{equation}
By Section \ref{par: case of two weights}, $\underline{\fu}(M/W_{p-1}M)$ is large if and only if this extension is totally nonsplit. It follows that we must have
\begin{equation}\label{eq27}
A \, + \, W_{p-1}M \ = \ W_pM.
\end{equation}
Indeed, otherwise, by modding out Eq. \eqref{eq26} by $A + W_{p-1}M$ we see that the pushforward of Eq. \eqref{eq25} to a nonzero subquotient of $Gr^W_pM$ splits, contradicting the fact that Eq. \eqref{eq25} is totally nonsplit (given by Hypothesis (ii)). 

Now consider the diagram
\[\begin{tikzcd}[remember picture]
& & A  \arrow[dr, ]  &   &  \\
0 \arrow[r] & W_{p-1}M \arrow[r, ] & W_{p}M \arrow[r, ] &  Gr^W_pM   \arrow[r] & 0 .
\end{tikzcd}
\begin{tikzpicture}[overlay,remember picture]
\path (\tikzcdmatrixname-1-3) to node[midway,sloped]{$\subset$}
(\tikzcdmatrixname-2-3);
\end{tikzpicture}
\]
We just saw that diagonal arrow is surjective. It follows that the extension in the diagram is the pushforward of an extension of $Gr^W_pM$ by $A\cap W_{p-1}M$ (under inclusion map). Thus the extension
\[
\mathcal{E}_{p-1}(W_pM) \ \in \ Ext^1(\mathbbm{1}, \inHom(Gr^W_pM,W_{p-1}M))
\]
is the pushforward of an extension of $\mathbbm{1}$ by 
\[
\inHom(Gr^W_pM,A\cap W_{p-1}M) \ \subset \  \inHom(Gr^W_pM,W_{p-1}M),
\]
i.e. that 
\[
\mathcal{E}_{p-1}(W_pM)/\inHom(Gr^W_pM,A\cap W_{p-1}M) 
\]
splits. By Theorem \ref{thm 1}, we get
\[
\underline{\fu}_{p-1}(W_pM) \ \subset \ \inHom(Gr^W_pM,A\cap W_{p-1}M).
\]
But since $\underline{\fu}(W_pM)$ is large, so is $\underline{\fu}_{p-1}(W_pM)$. Thus 
\[
\inHom(Gr^W_pM,A\cap W_{p-1}M) \ = \ \inHom(Gr^W_pM,W_{p-1}M).
\]
Since $Gr^W_pM$ is nonzero, this implies that $W_{p-1}M\subset A$. Combining with Eq. \eqref{eq27} we get that $A=W_pM$, as desired\footnote{Note that the assumption that $Gr^W_pM$ is nonzero is actually important for the proof. Thus when we want to apply Theorem \ref{thm3} to show that a given motive $M$ has a large $\underline{u}$, we do not have a choice about what to take as $p$; it is determined by the motive $M$.}.
\end{proof}

\begin{rem}
\begin{itemize}
\item[(1)] As pointed out in the proof, Hypothesis (ii) of Theorem \ref{thm3} is equivalent to the extension Eq. \eqref{eq25} being totally nonsplit. If we assume moreover that $Gr_pM$ is simple, then this is equivalent to $M/W_{p-1}M$ not being semisimple.
\item[(2)] Let $M$ be a motive which satisfies Eq. \eqref{eq28} (with $p<0$). It is easy to see that if $\underline{\fu}_p(M)$ is large, then so is $\underline{\fu}(M/W_{p-1}M)$. Indeed, if the latter is not large, then the pushforward of Eq. \eqref{eq25} to a nonzero quotient of $Gr^W_pM$ splits. The same split extension is then the pushforward of Eq. \eqref{eq26} to a nonzero quotient of $W_pM$, so that by Theorem \ref{thm 1} $\underline{\fu}_p(M)$ is not large. 

Now suppose that $\underline{\fu}(M)$ is large. As we observed in the beginning of the proof of Theorem \ref{thm3}, this implies that both $\underline{\fu}(W_pM)$ and $\underline{\fu}_p(M)$ are large. We record the conclusion:

\emph{If $M$ is  a motive satisfying Eq. \eqref{eq28} (with $p<0$) and $\underline{\fu}(M)$ is large, then both $\underline{\fu}(W_pM)$ and $\underline{\fu}(M/W_{p-1}M)$ are large. }   

Note that here we did not need to assume $M$ satisfies any independence axiom. Theorem \ref{thm3} asserts that if we further assume that $M$ satisfies the independence axiom given in Hypothesis (iii) of the theorem, then the converse to the statement above is also true.
 
\item[(3)] Hypothesis (iii) of the theorem (which was used in the proof to guarantee that $\mathcal{E}_p(M)/\underline{\fu}_p(M)$ splits) is actually important: the statement of the theorem is false if we remove Hypothesis (iii). See Subsection \ref{counterexample to Thm3 without IA hypothesis} for an example.
\end{itemize}
\end{rem}

\subsection{}\label{sec: compatible pairs, defn} In view of Theorem \ref{thm3} one may hope to form motives with large $\underline{\fu}$ by patching together suitable smaller such motives. The goal of the next few subsections is to try to classify, up to isomorphism, all motives $M$ with large $\underline{\fu}$ which satisfy Eq. \eqref{eq28} and which, up to isomorphism, have a fixed $W_{p-1}M$ (with large $\underline{\fu}$) and $Gr^W_pM$ (with the isomorphisms not part of the data). To this end, let us first consider a related  problem. For the discussion in this subsection, $\mathbf{T}$ can be any abelian category (we will eventually apply the discussion to our category of motives). 

Throughout, we fix objects $A, B$ and $C$ in $\mathbf{T}$ (in our final application, these will be respectively (the fixed objects which are to be isomorphic to) $Gr^W_pM$, $W_{p-1}M$, and $\mathbbm{1}$). Grothendieck considers the following problem in SGA 7 \cite[\S 9.3 of Expos\'{e} 9]{Gr68}: classify all tuples
\[
(M; (M_i)_{-3\leq i\leq 0}; \gamma_0, \gamma_{-1},\gamma_{-2}) 
\] 
where 
\[
M \, = \, M_0  \, \supset \,  M_{-1}  \, \supset \, M_{-2}  \, \supset \, M_{-3} \, = \,0
\]
are objects of $\mathbf{T}$ and 
\[
M/M_{-1} \, = \, M_0/M_{-1} \, \stackrel{\gamma_0}{\longrightarrow} \, C, \ \ M_{-1}/M_{-2} \, \stackrel{\gamma_{-1}}{\longrightarrow} \, A, \ \   M_{-2}/M_{-3} \, = M_{-2}\,  \stackrel{\gamma_{-2}}{\longrightarrow} \, B
\]
are isomorphisms. The classification is to be done up to isomorphisms of such tuples, defined in the obvious way. Here it is convenient for us to consider a slight variant of this problem, where we do not include the data of the isomorphisms $\gamma_i$ in the tuple, but instead just require that the quotients $M_0/M_{-1}$, $ M_{-1}/M_{-2}$ and $M_{-2}/M_{-3}=M_{-2}$ are isomorphic to $C$, $A$ and $B$, respectively. 

We say that a pair of extension classes
\[
(\mathcal{L},\mathcal{N}) \ \in \  Ext^1(A, B) \, \times \, Ext^1(C, A)
\]
is {\it compatible} if there is a commutative diagram in $\mathbf{T}$:
\begin{equation}\label{eq18}
\begin{tikzcd}
   & & 0 \arrow{d} & 0 \arrow{d} &\\
   0 \arrow[r] & B \ar[equal]{d} \arrow[r, ] & L  \arrow[d, ] \arrow[r, ] &  A \arrow{d} \arrow[r] & 0 \\
   0 \arrow[r] & B \arrow[r, ] & M \arrow{d} \arrow[r, ] &  N \arrow{d}  \arrow[r] & 0 \, ,\\
   & & C \arrow{d} \ar[equal]{r} & C \arrow{d} & \\
   & & 0 & 0 &   
\end{tikzcd}
\end{equation}
where the rows and columns are exact, the first (complete) row represents $\mathcal{L}$, and the second (complete) column represents $\mathcal{N}$. We say an object $M$ is {\it attached to} the pair $(\mathcal{L},\mathcal{N})$ if it fits in a diagram as above. Note that if we have a diagram as above, (by adjusting the maps when needed) we may replace the first row (resp. second column) by any other representative of $\mathcal{L}$ (resp. $\mathcal{N}$). 

In SGA 7, a diagram as above is called an {\it extension panach\'{e}e} \footnote{or as Bertrand translates in \cite{Ber98}, a {\it blended} extension} of the second column sequence by the top row sequence. Thus to say the pair $(\mathcal{L},\mathcal{N})$ is compatible amounts to saying that an {\it extension panach\'{e}e} of (an or every representative of) $\mathcal{N}$ by (an or every representative of) $\mathcal{L}$ exists, or that $(\mathcal{L},\mathcal{N})$ is ``{\it panachable}", in the language of \cite{Ber13}.

The theory of Yoneda extensions gives a simple characterization of compatible pairs. Let 
\[
\circ \,: \,  Ext^1(A, B) \, \times \, Ext^1(C, A)  \ \longrightarrow \ Ext^2(C, B) 
\]
be the Yoneda (composition) pairing; it sends the pair $(\mathcal{L}, \mathcal{N})$ with $\mathcal{L}$ given by
\begin{equation}\label{eq17}
0 \ \longrightarrow \  B  \ \longrightarrow \ L \ \stackrel{\pi}{\longrightarrow} \ A  \ \longrightarrow \ 0
\end{equation}
and $\mathcal{N}$ given by
\[
0 \ \longrightarrow \  A  \ \stackrel{\iota}{\longrightarrow} \ N \ \longrightarrow \ C  \ \longrightarrow \ 0
\]
to the extension $\mathcal{L}\circ \mathcal{N}$ given by
\[
0 \ \longrightarrow \  B  \ \longrightarrow \ L \ \stackrel{\iota\circ \pi}{\longrightarrow} \ N  \ \longrightarrow \ C  \ \longrightarrow \ 0.
\]

\begin{lemma}\label{criterion for compatibility of a pair} \ \\
\vspace{-.1in}
\begin{itemize}
\item[(a)] The pair $(\mathcal{L}, \mathcal{N})$ is compatible if and only if $\mathcal{L}\circ \mathcal{N} \ = \ 0$.
\item[(b)] Suppose $Ext^1(C,B)=0$. If $(\mathcal{L}, \mathcal{N})$ is compatible, then up to isomorphism there is a unique object attached to it.
\end{itemize}
\end{lemma}

\begin{proof}
This is Lemma 9.3.8 of \cite{Gr68}. Fix the extension Eq. \eqref{eq17} representing $\mathcal{L}$. If $M$ is attached to the pair, fitting into a diagram as in Eq. \eqref{eq18}, then the class 
\[\mathcal{M} \ \in \ Ext^1(C, L)\]
of the first column in the diagram pushes forward to $\mathcal{N}$ under $\pi$. Conversely, if $\mathcal{N}$ is in the image of the pushforward
\[
\pi_\ast \, : \, Ext^1(C,L)  \ \longrightarrow \ Ext^1(C,A),
\]
with $\mathcal{M}$ represented by
\[
0   \ \longrightarrow \ L   \ \longrightarrow \  M   \ \longrightarrow \ C  \ \longrightarrow \ 0
\]
in the preimage of $\mathcal{N}$, then the object $M$ is attached to our pair. Thus the pair $(\mathcal{L}, \mathcal{N})$ is compatible if and only if $\mathcal{N}$ is in the image of $\pi_\ast$. Now by the general theory of Yoneda extensions, applying the functor $Hom(C, -)$ to  Eq. \eqref{eq17} we get an exact sequence
\[
 Ext^1(C,B)  \ \longrightarrow \  Ext^1(C,L)  \ \stackrel{\pi_\ast}{\longrightarrow} \ Ext^1(C,A)  \ \xrightarrow{\delta = \mathcal{L} \circ -} \ Ext^2(C,B)  
\]
(see \cite[Section 3]{Bu59} or \cite[page 561]{Yo60}). This proves Part (a).

As for the statement in Part (b), if $M$ and $M'$ are attached to $(\mathcal{L},\mathcal{N})$, fitting into diagrams as in Eq. \eqref{eq18} with the classes of the corresponding first columns denoted by $\mathcal{M}$ and $\mathcal{M}'$ (both in $Ext^1(C,L)$) respectively, then it follows from the above long exact sequence that $\mathcal{M}$ and $\mathcal{M}'$ differ by an element in the image of $Ext^1(C,B)$. If this Ext group is zero, then $\mathcal{M}=\mathcal{M}'$, and hence in particular $M$ and $M'$ are (non-canonically) isomorphic.
\end{proof}

\subsection{}\label{sec: compatible pairs II} We shall continue in the setting of the previous subsection ($\mathbf{T}$ any abelian category, and $B, A, C$ three fixed objects of $\mathbf{T}$). Our goal in this subsection is to see when the same object is attached to two compatible pairs of extensions. 

We use the notation $End( \ )$ (resp. $Aut( \ )$) for the endomorphism algebra (resp. automorphism group) of an object in $\mathbf{T}$. The endomorphism algebra $End(A)$ of $A$ acts on both $Ext^1(A, B)$ and $Ext^1(C,A)$. Indeed, the action on $Ext^1(A, B)$ is a right action given by pullback: if $f$ is an endomorphism of $A$, set $\mathcal{L}\cdot f:=f^\ast\mathcal{L}$ ($f^\ast$ for pullback along $f$). The action on $Ext^1(C,A)$ is a left action given by push forward: $f\cdot \mathcal{N} := f_\ast\mathcal{N}$ (to see the bilinearity properties of these actions see \cite{Bu59} or \cite{Yo60}). If $f$ is an automorphism of $A$, then $\mathcal{L}\cdot f$ and $f\cdot \mathcal{N}$ are simply obtained by twisting respectively the surjective and injective arrows of $\mathcal{L}$ and $\mathcal{N}$ by $f^{-1}$, i.e. $\mathcal{L}\cdot f$ (resp. $f\cdot \mathcal{N}$) is the class of the extension obtained by replacing the surjective (resp. injective) arrow $\pi$ (resp. $\iota$) in a representative of $\mathcal{L}$ (resp. $\mathcal{N}$) by $f^{-1}\circ \pi$ (resp. $\iota\circ f^{-1}$). 

We restrict the two actions above on $Ext^1(A, B)$ and $Ext^1(C,A)$ to the actions of the group $Aut(A)$. We also modify the action on $Ext^1(C,A)$ so that it also becomes a right action, by setting $\mathcal{N}\cdot f := f^{-1}_\ast\mathcal{N}$. Thus $\mathcal{N}\cdot f$ is the class of the extension obtained by twisting the injective arrow of $\mathcal{N}$ by $f$.  Similarly, we have right actions of $Aut(B)$ (resp. $Aut(C)$) on $Ext^1(A, B)$ (resp. $Ext^1(C,A)$). 

We now equip the product 
\begin{equation}\label{eq21}
Ext^1(A, B) \, \times \, Ext^1(C, A)
\end{equation}
with the following right actions of $Aut(B)$, $Aut(A)$, and $Aut(C)$: the group $Aut(B)$ (resp. $Aut(C)$) acts by acting on the first (resp. second) factor, and $Aut(A)$ acts diagonally, i.e. by the formula
\[
(\mathcal{L},\mathcal{N})\cdot f \ := \ (\mathcal{L}\cdot f,\mathcal{N}\cdot f) \ = \ (f^\ast\mathcal{L},f^{-1}_\ast\mathcal{N}).
\] 
There three actions commute with one another. Indeed, the actions of $Aut(B)$ and $Aut(C)$ trivially commute, and the commutativity of the actions of $Aut(A)$ with each of $Aut(B)$ and $Aut(C)$ is clear from the description of the actions in terms of twisting the arrows, as different groups act by twisting different arrows. Thus we get an action of $Aut(B)\times Aut(A) \times Aut(C)$ on the product Eq. \eqref{eq21}. We say two pairs of extensions are {\it equivalent} if they are in the same orbit of this action.

\begin{lemma}\label{equivalent pairs} Let $(\mathcal{L},\mathcal{N})$ and $(\mathcal{L}',\mathcal{N}')$ be in Eq. \eqref{eq21}. 
\begin{itemize}
\item[(a)] Suppose $(\mathcal{L},\mathcal{N})$ and $(\mathcal{L}',\mathcal{N}')$ are equivalent. Then every object attached to the pair $(\mathcal{L},\mathcal{N})$ is also attached to the pair $(\mathcal{L}',\mathcal{N}')$. (In particular, $(\mathcal{L},\mathcal{N})$ is compatible if and only if $(\mathcal{L}',\mathcal{N}')$ is compatible.)
\item[(b)] Suppose every object of $\mathbf{T}$ is equipped with an exact functorial increasing filtration $W_\cdot$ which is finite on every object (we refer to this as the weight filtration). Suppose moreover that the highest weight of $B$ is less than the lowest weight of $A$, and that the highest weight of $A$ is less than the lowest weight of $C$. Then if there is an object $M$ attached to both $(\mathcal{L},\mathcal{N})$ and $(\mathcal{L}',\mathcal{N}')$, then the two pairs are equivalent.
\end{itemize}
\end{lemma}

\begin{proof} (a) Let $(\mathcal{L}',\mathcal{N}')=(\mathcal{L},\mathcal{N})\cdot (f_B, f_A, f_C)$ for some $f_B\in Aut(B)$, $f_A\in Aut(A)$, and $f_C\in Aut(C)$. Suppose $M$ is attached to $(\mathcal{L},\mathcal{N})$. In a diagram as in Eq. \eqref{eq18} (with the first row and second column respectively representing $\mathcal{L}$ and $\mathcal{N}$), twist the arrows $B\longrightarrow L$ and $B\longrightarrow M$ by $f_B$, the arrow $L\longrightarrow A$ by $f_A^{-1}$, the arrow $A\longrightarrow N$ by $f_A$, and the arrows $M\longrightarrow C$ and $N\longrightarrow C$ by $f_C$, while keeping $L\longrightarrow M$ and $M\longrightarrow N$ unchanged. The diagram remains commutative and with exact rows and columns, and its first row (resp. second column) represents $\mathcal{L}'$ (resp. $\mathcal{N}'$).

(b) Suppose an object $M$ is attached to both $(\mathcal{L},\mathcal{N})$ and $(\mathcal{L}',\mathcal{N}')$. We consider two diagrams as in Eq. \eqref{eq18}, one with objects $L, N$ with the first row and second column representing $\mathcal{L}$ and $\mathcal{N}$, and the other with objects $L', N'$ with the first row and second column representing $\mathcal{L}'$ and $\mathcal{N}'$. In the diagram for $(\mathcal{L},\mathcal{N})$, we name the maps as follows: In the first row, (resp. second row, second column) the injective arrow is $\iota_L$ (resp. $\iota_M$, $\iota_N$) and the surjective arrow is $\pi_L$ (resp. $\pi_M$, $\pi_N$). We refer to the maps $L\longrightarrow M$ and $M\longrightarrow C$ as $\alpha$ and $\beta$, respectively. Accordingly, denote the maps in the diagram for $(\mathcal{L}',\mathcal{N}')$ by $\iota_{L'},\pi_{L'}, \iota'_M, \pi_M', \iota_{N'}, \pi_{N'}, \alpha'$ and $\beta'$ (each map being the analogue to its lookalike in the first diagram). (Note that the central object in both diagrams in $M$.)

Let $b, a$ and $c$ be respectively the highest weights of $B, A$ and $C$. Focusing on the first diagram, using exactness of the weight filtration together with the hypothesis that every weight of $B$ is less than every weight of $A$, which in turn is less than every weight of $C$, we see that
\[
W_bL \ = \ \iota_L(B), \quad W_aL \ = \ L, \quad W_bN \ = \ 0, \quad W_aN \ = \ \iota_N(A) ,\quad W_cN \ = \ N
\]
and 
\[
W_bM \ = \ \iota_M(B), \quad W_aM \ = \ \alpha(L), \quad W_cM \ = \ M. 
\]
We have similar equalities for the $'$-adorned analogues coming from the second diagram. In particular,
\[
\iota_M(B) \ = \ \iota'_M(B) \ = \ W_bM, \quad\quad \alpha(L) \ = \ \alpha'(L') \ = \ W_aM.
\]
Thus we get isomorphisms $\alpha^{-1}\alpha': L' \longrightarrow L$ and ${\iota}^{-1}_M\iota'_M : B\longrightarrow B$ (uniquely) defined by the property that ${\alpha}(\alpha^{-1}\alpha')=\alpha'$ and $\iota_M({\iota}^{-1}_M\iota'_M)=\iota'_M$. We have a commutative diagram with exact rows
\[
\begin{tikzcd}[column sep = large, row sep = large ]
 0 \arrow[r] & B \arrow{d}{{\iota}^{-1}_M\iota'_M} \arrow{r}{ \iota_{L'}} & L'  \arrow{d}{ \alpha^{-1}\alpha'} \arrow{r}{ \pi_{L'}} &  A \arrow{d}{=: \, \gamma} \arrow[r] & 0 \\
 0 \arrow[r] & B \arrow{r}{\iota_{L}} & L  \arrow{r}{\pi_{L}} &  A   \arrow[r] & 0 \, , 
\end{tikzcd}
\]
where the vertical arrows are isomorphisms (to see the commutativity of the first square further compose with $\alpha$). Thus $\mathcal{L}'$ is obtained from $\mathcal{L}$ by twisting $\iota_L$ by ${\iota}^{-1}_M\iota'_M$ and twisting $\pi_L$ by $\gamma^{-1}$.

On the other hand, since we have $\iota_M(B) = \iota'_M(B) = W_bM$, by exactness of the second rows in the diagrams of the two pairs, $\pi_M$ and $\pi'_M$ induce isomorphisms
\[
\overline{\pi_M} \, : \, M/W_bM \ \longrightarrow \ N \, ,  \quad \quad \quad \overline{\pi'_M} \, : \, M/W_bM \ \longrightarrow \ N'.
\]
Similarly, thanks to exactness of the first columns (and on recalling $\alpha(L) = \alpha'(L')  =  W_aM$), we have isomorphisms
\[
\overline{\beta} \, : \, M/W_aM \ \longrightarrow \ C \, , \quad \quad \quad  \overline{\beta'} \, : \, M/W_aM \ \longrightarrow \ C,
\]
induced by $\beta$ and $\beta'$, respectively. We now have a commutative diagram with exact rows
\begin{equation}\label{eq22}
\begin{tikzcd}[column sep = large, row sep = large ]
 0 \arrow[r] & A \arrow{d}{=: \, \lambda  } \arrow{r}{\iota_{N'} } & N'  \arrow{d}{ \overline{\pi_M}{\overline{\pi'_M}}^{-1}} \arrow{r}{\pi_{N'}} &  C \arrow{d}{\overline{\beta}\, {\overline{\beta'}}^{-1} } \arrow[r] & 0 \\
 0 \arrow[r] & A \arrow{r}{\iota_{N}} & N  \arrow{r}{\pi_{N}} &  C   \arrow[r] & 0 ,
\end{tikzcd}
\end{equation} 
where the vertical arrows are isomorphisms (to see commutativity of the second square pre-compose with $\pi'_M: M\longrightarrow N'$). It follows that $\mathcal{N'}$ is obtained from $\mathcal{N}$ by twisting $\iota_N$ by $\lambda$ and twisting $\pi_N$ by $\overline{\beta'}\,{\overline{\beta}}^{-1}$. 

To complete the proof, it suffices to show that $\gamma=\lambda$, as then
\[
(\mathcal{L}',\mathcal{N}') \ = \ (\mathcal{L},\mathcal{N})\cdot({\iota}^{-1}_M\iota'_M, \gamma,\overline{\beta}\, {\overline{\beta'}}^{-1} ).
\]
Ignoring the dashed arrow, we have a commutative diagram
\[
\begin{tikzcd}
C & & \arrow[ll, "\overline{\beta}" ' , "\simeq"  ] M/W_aM \arrow[rr, "\overline{\beta'}" , "\simeq" '] & & C \\
N \arrow[u, "\pi_{N}"]  & & \arrow[ll, "\overline{\pi_M}" ' , "\simeq"  ] M/W_bM \arrow[u] \arrow[rr, "\overline{\pi'_M}" , "\simeq" '] & & N' \arrow[u, "\pi_{N'}"] \\  
  & & M \arrow [ull, "\pi_M" ] \arrow[u] \arrow[urr, "\pi'_M" ' ] & & \\  
    & L \arrow[dl, "\pi_L" ' ] \arrow[ur, "\alpha" ] & & \arrow[ll, " \alpha^{-1}\alpha' " ' , "\simeq" ]   \arrow[ul, "\alpha' " '] L' \arrow[dr, " \pi_{L'}" ] & \\  
  A \arrow[uuu, " \iota_N " ] & & & & A  \arrow[llll, dashed] \arrow[uuu, " \iota_{N'} " ] \, , 
 \end{tikzcd}
\]
where the vertical arrows in the middle are the obvious maps. The map $\gamma$ is the unique map which if it is included as the dashed arrow, it makes the bottom trapezoid of the diagram commute. But from the diagram we easily see that $\lambda$ also does this job. Indeed, to check commutativity of the trapezoid with $\lambda$ as the dashed arrow, it is enough to check commutativity after composing with $\iota_N$. Now using commutativity of the rest of the diagram above and the left square in Eq. \eqref{eq22}, we have 
\[
\iota_N \, \pi_L \, (\alpha^{-1}\alpha') \ = \ \pi_M \, \alpha' \ = \ \iota_N \, \lambda \, \pi_{L'} \,.
\]
\end{proof}

\subsection{}\label{sec: classification result compatible pairs vs objects with ... up to isomorphism} We now combine the results of the previous two subsections on compatible pairs. We shall assume that $\mathbf{T}$  is an abelian category equipped with a weight filtration (i.e. a functorial, exact, increasing filtration which is finite on every object). As in the previous two subsections, $B, A, C$ are fixed objects of $\mathbf{T}$. The following result, which for future reference we record as a proposition, has been mostly already proved in the previous two subsections.

\begin{prop}\label{prop compatible pairs}
Suppose every weight of $B$ is less than every weight of $A$, and that every weight of $A$ is less than every weight of $C$. Let $b,a,c$ be the highest weights of $B, A, C$, respectively. 
\begin{itemize}
\item[(a)] Any pair of extensions $(\mathcal{L},\mathcal{N})$ in Eq. \eqref{eq21} is compatible if and only if $\mathcal{L}\circ \mathcal{N} = 0$ in $Ext^2(C,B)$.
\item[(b)] If $M$ is an object which is attached to some pair of extensions in Eq. \eqref{eq21}, then we have
\begin{equation}\label{eq23}
B \ \simeq \ W_bM ,\, \quad  A  \ \simeq \ W_aM/W_bM, \, \quad C  \ \simeq \ M/W_aM.
\end{equation}
\item[(c)] Any object $M$ satisfying Eq. \eqref{eq23} is attached to some pair $(\mathcal{L},\mathcal{N})$ of extensions in Eq. \eqref{eq21}. Moreover, $M$ is attached to any other pair $(\mathcal{L}',\mathcal{N}')$ if and only if $(\mathcal{L}',\mathcal{N}')$ is equivalent to $(\mathcal{L},\mathcal{N})$. We have a (well-defined) surjective map
\begin{align*}
\text{the collection of objects $M$ satisfying}   \ \quad &  \ \ \raise-1.5ex\hbox{$\longrightarrow$} & \text{the collection of compatible pairs} \\
\text{ Eq. \eqref{eq23}, up to isomorphism}  \ \quad & &   \text{ in Eq. \eqref{eq21}, up to equivalence}
\end{align*}
which sends the isomorphism class of $M$ to the equivalence class of any pair (or all pairs) $(\mathcal{L},\mathcal{N})$ to which $M$ is attached.
\item[(d)] If $Ext^1(C, B)=0$, then the surjection above is a bijection.
\end{itemize}
\end{prop}

\begin{proof}
(a) This is Lemma \ref{criterion for compatibility of a pair}(a).

(b) This follows from the observations made at the beginning of the proof of Lemma \ref{equivalent pairs}(b) about the weight filtration of $M$. (Note that the isomorphisms are non-canonical, as they depend on the particular choice of diagram Eq. \eqref{eq18}.)

(c)  Given $M$ satisfying Eq. \eqref{eq23}, we have a diagram
\[
\begin{tikzcd}
   & & 0 \arrow{d} & 0 \arrow{d} &\\
   0 \arrow[r] & W_{b}M \ar[equal]{d} \arrow[r, ] & W_aM  \arrow[d, ] \arrow[r, ] &  W_aM/W_bM \arrow{d} \arrow[r] & 0 \\
   0 \arrow[r] & W_{b}M \arrow[r, ] & M \arrow{d} \arrow[r, ] &  M/W_{b}M  \arrow{d}  \arrow[r] & 0 \\
   & & M/W_aM \arrow{d} \ar[equal]{r} & M/W_aM \arrow{d} & \\
   & & 0 & 0 &   
\end{tikzcd}
\]
(with obvious maps, exact rows and columns). Now use isomorphisms as in Eq. \eqref{eq23} to replace $W_{b}M$, $W_aM/W_bM$, and $M/W_aM$ respectively by $B$, $A$, and $C$. Take $\mathcal{L}$ (resp. $\mathcal{N}$) to be the extension class of the top row (resp. last column) in the new diagram. Then $M$ is attached to the (compatible) pair $(\mathcal{L},\mathcal{N})$. By Lemma \ref{equivalent pairs}(b), $M$ is attached to another pair $(\mathcal{L}',\mathcal{N}')$ if and only if $(\mathcal{L}',\mathcal{N}')$ is equivalent to $(\mathcal{L},\mathcal{N})$. On the other hand, if $M'$ is isomorphic to $M$, then $M'$ is clearly attached to the same pairs as $M$. Thus we have a well-defined map as in the statement. It is surjective by the definition of compatibility and Part (b).

(d) This follows from Lemma \ref{criterion for compatibility of a pair}(b) and Lemma \ref{equivalent pairs}(a).
\end{proof}

\subsection{}\label{sec: summary of results of the third part} We now return to the discussion of motives with large $\underline{\fu}$ (with $\mathbf{T}$ again a Tannakian category of mixed motives or the category of rational mixed Hodge structures). Given any two motives $A$ and $B$, let us say an extension class in $Ext^1(A,B)$ has a large $\underline{\fu}$ if the object in the middle of a representing short exact sequence has a large $\underline{\fu}$. This is clearly well-defined, and moreover, the property of having a large $\underline{\fu}$ is invariant under the action of $Aut(A)\times Aut(B)$ (because the collection of the objects that can appear as the middle object for two extension classes in the same orbit are the same, as by twisting the arrows we can turn a representative of one extension class to a representative of another extension class in the same orbit). Note that if $A$ is simple (resp. pure), then an extension class in $Ext^1(\mathbbm{1},A)$ has a large $\underline{\fu}$ if and only if it is nonsplit (resp. totally nonsplit). 

We say a pair of extensions $(\mathcal{L},\mathcal{N})$ in Eq. \eqref{eq21} has a large $\underline{\fu}$ if both extensions $\mathcal{L}$ and $\mathcal{N}$ have a large $\underline{\fu}$. This property is invariant under our notion of equivalence of pairs.

We now fix an integer $p<0$, and motives $B$ and $A$ with 
\[
B \ = \ W_{p-1}B , \quad A\cong Gr^W_pA \ \neq \ 0.
\]
(In other words, all weights of $B$ are $<p$, and $A$ is nonzero and pure of weight $p$; note that $B$ may be mixed.) Proposition \ref{prop compatible pairs} gives a surjection (bijection if $Ext^1(\mathbbm{1},B)=0$) from the collection of motives $M$ satisfying 
\begin{equation}\label{eq105}
W_{p-1}M \ \simeq \ B, \quad Gr_p^WM \ \simeq A, \quad M/W_pM \ \simeq \ \mathbbm{1}
\end{equation}
up to isomorphism to the collection of compatible pairs in 
\[
Ext^1(A,B) \, \times \, Ext^1(\mathbbm{1},A) 
\]
(= the kernel of the composition pairing into $Ext^2(\mathbbm{1},B)$) up to equivalence (i.e. the action of $Aut(B)\times Aut(A) \times Aut(\mathbbm{1})$). By Theorem \ref{thm3}, if $B\oplus A \oplus \mathbbm{1}$ satisfies the independence axiom $(IA1)_{\{p,q\}}$ for every $q\leq p$, then given any compatible pair $(\mathcal{L},\mathcal{N})$ with a large $\underline{\fu}$, any object $M$ attached to the pair also has a large $\underline{\fu}$. Conversely, if an object $M$ satisfying Eq. \eqref{eq105} has a large $\underline{\fu}$, then so does any pair $(\mathcal{L},\mathcal{N})$ in the equivalence class of the extension pairs corresponding to $M$ (see Remark (2) after Theorem \ref{thm3}; note that here no independence axiom needs to be satisfied).

We record the following special case as a corollary:

\begin{cor}\label{cor: motives with three weights and Gr^W isomorphic to Q(n)+A+1, bijective case }
Let $-2n<p<0$ and $p\neq -n$. Let $A$ be a nonzero simple motive of weight $p$. Suppose moreover that $Ext^1(\mathbbm{1},\QQ(n))=0$. Then there is a bijection 
\begin{align*}
\text{the collection of objects $M$ \ \ }   \ \quad &   & \text{the collection of compatible pairs} \\
\text{with $Gr^WM\, \simeq  \, \QQ(n)\oplus A\oplus \mathbbm{1}$ \quad }   &   \longrightarrow &   \text{of nonsplit extentions in \quad \quad \quad \ }\\
\text{and large $\underline{\fu}(M)$, up to  \quad \quad \quad \, \,}   & &   \text{$Ext^1(A, \QQ(n)) \, \times \, Ext^1(\mathbbm{1}, A)$, \ }\\ 
 \text{isomorphism \quad \quad \quad \quad \quad \quad \quad \ \ }  & &   \text{up to equivalence, \quad\quad \quad \quad \quad \quad }
\end{align*}
which assigns to the isomorphism class of an object $M$ the equivalence class of the compatible pairs to which $M$ is attached. If we omit the condition $Ext^1(\mathbbm{1},\QQ(n))=0$, this map is well-defined and surjective.
\end{cor}

(Note that the condition $p\neq -n$ guarantees $(IA3)$.)

\begin{rem} \ \\
\vspace{-.1in}
\begin{itemize}
\item[(1)] Recall that in the category $\mathbf{MHS}$ of rational mixed Hodge structures, the $Ext^2$ (and hence, all higher Ext) groups vanish (see \cite{Be83}). The $Ext^1$ groups in this case are described by the results of Carlson \cite{Ca80}.

\item[(2)] Let $\mathbf{MT}(K)$ be Voevodsky's category of mixed Tate motives over a number field $K$. The $Ext^2$ groups in $\mathbf{MT}(K)$ are zero, and the groups
\[
Ext^1_{\mathbf{MT}(K)}(\mathbbm{1}, \QQ(n))
\]
are given by the K-theory of the field $K$ modulo torsion, which in turn is described by theorems of Borel and Soul\'{e} (and Dirichlet in the case of $K_1$). In particular, if $K$ is totally real and $n$ is even, the $Ext^1$ group above vanishes. (See \cite{DG05} for the precise description of the Ext groups in $\mathbf{MT}(K)$ and the subcategory of mixed Tate motives over the ring of integers of $K$. Note that if $\mathbf{MM}(K)$ is any category of mixed motives over $K$ for which the full Tannakian subcategory generated by $\QQ(1)$ and closed under subobjects and extensions is equivalent to Voevodsky's $\mathbf{MT}(K)$, then the $Ext^1$ groups above are the same in $\mathbf{MM}(K)$ and $\mathbf{MT}(K)$.)

\item[(3)] In a good category of mixed motives over a number field, the $Ext^2$ groups are expected to vanish. The $Ext^1$ groups in such a category should be related to Chow groups and motivic cohomology (and algebraic K-theory). See for instance, Nekovar's expository article \cite{Ne94} or Jannsen's \cite{Ja94}. 

\end{itemize}
\end{rem}

\subsection{}\label{sec: classification of 3 dimensional mixed Tate motives with large u} In this subsection, we shall take $\mathbf{T}$ to be Voevodsky's category $\mathbf{MT}(\QQ)$ of mixed Tate motives over $\QQ$. As an application of the previous results, we shall classify (up to isomorphism) all 3-dimensional objects of $\mathbf{MT}(\QQ)$ with three distinct weights, large $\underline{\fu}$, and satisfying an independence axiom (see below for more details).\footnote{The classification is then valid in any Tannakian category $\mathbf{MM}(\QQ)$ of mixed motives over $\QQ$ for which the smallest full Tannakian subcategory containing $\QQ(1)$ and closed under subobjects and extensions is equivalent to $\mathbf{MT}(\QQ)$.} Note that for any 3-dimensional object $M$ of $\mathbf{MT}(\QQ)$ with three distinct weights and large $\underline{\fu}$, the unipotent radical of the motivic Galois group $\mathcal{G}(M,\omega_B)$ (with $\omega_B$ the Betti realization functor) has dimension equal to 3 ( = $\dim W_{-1}End(\omega_BM)$). Since
\[
\mathcal{G}(Gr^WM,\omega_B) \  \simeq \ \mathbb{G}_m,
\]
the motivic Galois group $\mathcal{G}(M,\omega_B)$ has dimension 4. Thus Grothendieck's period conjecture would predict that the transcendence degree of the field generated by the periods of $M$ should be 4. 

Let us first recall the description of the Ext groups between simple objects in $\mathbf{MT}(\QQ)$ (see \cite{DG05}, for instance):
\[
\dim Ext^1(\mathbbm{1},  \QQ(n)) \ = \ \begin{cases}
\ 0 \hspace{.3in} &\text{if $n$ is even or $\leq 0$}\\
\ 1 &\text{if $n$ is odd and $\geq 3$}
\end{cases}
\]
\begin{equation}\label{eq19}
Ext^1(\mathbbm{1},  \QQ(1)) \ \cong \ \QQ^\times \otimes \QQ \, .
\end{equation}
Moreover, $Ext^2$ groups all vanish in $\mathbf{MT}(\QQ)$. 

Back to our classification problem, we may assume that our motives have highest weight equal to 0. We shall classify all motives with an associated graded of the form
\[
\QQ(n) \oplus \QQ(k) \oplus \mathbbm{1} \hspace{.5in} (n>k>0, \ n\neq 2k)
\]
which have large $\underline{\fu}$. (The condition $n\neq 2k$ is an independence axiom. The case where $n=k$ is complicated, as then one can no longer use Theorem \ref{thm3}.) For any such motive, the pair $(\mathcal{L},\mathcal{N})$ in
\begin{equation}\label{eq29}
Ext^1(\QQ(k), \QQ(n)) \, \times \, Ext^1(\mathbbm{1}, \QQ(k))
\end{equation}
associated to it by Corollary \ref{cor: motives with three weights and Gr^W isomorphic to Q(n)+A+1, bijective case } (also see Proposition \ref{prop compatible pairs}) has nonsplit entries. In view of the description of the $Ext^1$ groups in the category, we see that $k$ must be odd and $n$ must be even. But then we have a bijection as in Corollary \ref{cor: motives with three weights and Gr^W isomorphic to Q(n)+A+1, bijective case }. 

Let us consider the action of $Aut(\QQ(n))\times Aut(\QQ(k)) \times Aut(\mathbbm{1})$ on Eq. \eqref{eq29}. Since the automorphism group of every $\QQ(a)$ is $\QQ^\ast$, it follows from bilinearity of the actions of $End(A)$ on $Ext^1(A,B)$ and $Ext^1(B,A)$ (for any $A,B$ in any $K$-linear category) that the action of $Aut(\QQ(k))$ can be absorbed into the actions of the other two factors: $(\lambda,\gamma,\delta)$ acts the same as $(\lambda\gamma^{-1},1,\gamma\delta)$ (where $\lambda,\gamma,\delta\in\QQ^\ast$). It follows that an orbit of the action of $Aut(\QQ(n))\times Aut(\QQ(k)) \times Aut(\mathbbm{1})$ of Eq. \eqref{eq29} coincides with an element of  
\begin{equation}\label{eq30}
\bigm(Ext^1(\QQ(k), \QQ(n))\bigm/Aut(\QQ(n))\bigm) \  \times \ \bigm(Ext^1(\mathbbm{1}, \QQ(k))\bigm/  Aut(\mathbbm{1})\bigm)
\end{equation}
(with both actions made right actions, as before).

\underline{Case I}: $k=1$. Then $n$ is $\geq 4$ (and even), and 
\[
Ext^1(\QQ(k),\QQ(n)) \ \cong \ Ext^1(\mathbbm{1},\QQ(n-k))
\]
is a 1-dimensional vector space over $\QQ$, and all its nonzero elements are in the same $Aut(\QQ(n))$-orbit. 

The extensions of $\mathbbm{1}$ by $\QQ(1)$ are the Kummer motives. For each positive rational number $r$, let 
\[
[r] \ \in Ext^1(\mathbbm{1},  \QQ(1)) 
\]
be the extension class arising from the weight filtration of the 1-motive (see \cite{De74}) 
\begin{equation}\label{eq31}
K_r \ := \ [\ZZ\, \xrightarrow{1\mapsto r} \, \mathbb{G}_m]
\end{equation}
(considered as an object of $\mathbf{MT}(\QQ)$). Then $[r]$ is the element of $Ext^1(\mathbbm{1},\QQ(1))$ corresponding to $r\otimes 1$ under the isomorphism Eq. \eqref{eq19}. Thus $\{[p]: \text{$p$ prime $>0$}\}$ is a basis of $Ext^1(\mathbbm{1},\QQ(1))$ (over $\QQ$). A complete inequivalent set of representatives for the nonzero orbits of the action of $\QQ^\ast=Aut(\mathbbm{1})$ on $Ext^1(\mathbbm{1},  \QQ(1))$ is formed by the elements $[r]$, where $r$ runs through all rational numbers $>1$ which are not of the form $s^a$ for any $s\in\QQ$ and $a\in\ZZ$ with $a>1$. In view of Corollary \ref{cor: motives with three weights and Gr^W isomorphic to Q(n)+A+1, bijective case }, each such $[r]$ gives a (unique, up to isomorphism) motive $M_{n,r}$ with large $\underline{\fu}$ and associated graded isomorphic to
\[
\QQ(n) \oplus \QQ(1) \oplus \mathbbm{1}.
\]
These motives are non-isomorphic, and are up to isomorphism, all the motives with associated graded as above and large $\underline{\fu}$. 

A discussion of the periods of $M_{n,r}$ is in order. By construction, $W_{-2}M_{n,r}$ is a nontrivial extension of $\QQ(1)$ by $\QQ(n)$. Being a twist (by $\QQ(1)$) of a nontrivial extension of $\mathbbm{1}$ by $\QQ(n-1)$, the motive $W_{-2}M_{n,r}$ has the period matrix
\[
\begin{pmatrix}
(2\pi i)^{-n} & (2\pi i)^{-n}\zeta(n-1)  \\
0 & (2\pi i)^{-1}\end{pmatrix}
\]
with respect to suitably chosen bases of Betti and de Rham realizations. (Note that $n-1$ is odd and $\geq 3$. That a nontrivial extension of $\mathbbm{1}$ by $\QQ(n-1)$ has $\zeta(n-1)/(2\pi i)^{n-1}$ as a period follows from the work \cite{De89} of Deligne in the setting of realizations, and later the work \cite{DG05} of Deligne and Goncharov in the setting of Voevodsky motives.) One the other hand, $M_{n,r}$ has the Kummer 1-motive $K_r$ as a subquotient (by $W_{-2n}M_{n,r}=W_{-3}M_{n,r}$). With respect to suitably chosen bases of Betti and de Rham realizations, $K_r$ has the period matrix
\[
\begin{pmatrix}
(2\pi i)^{-1} & (2\pi i)^{-1}\log r  \\
0 & 1\end{pmatrix}
\]
(see \cite{De74} for the explicit realizations of 1-motives). With respect to suitably chosen bases, the matrix of periods of $M_{n,r}$ looks like
\[
\begin{pmatrix}
(2\pi i)^{-n}&(2\pi i)^{-n}\zeta(n-1) &\ast\\
0&(2\pi i)^{-1}&(2\pi i)^{-1}\log r \\
0&0&1
\end{pmatrix}.
\]
As mentioned earlier, Grothendieck's period conjecture predicts the transcendence degree of the field generated over $\QQ$ by the periods of $M_{n,r}$ to be 4. Thus assuming the period conjecture, the numbers 
\[
2\pi i, \, \log r, \, \zeta(n-1), \, \text{and the entry denoted by $\ast$}
\]
are algebraically independent over $\QQ$. 

It would be very interesting to somehow calculate the entry $\ast$ in the matrix above. As we discussed in the Introduction, when $r\neq 2$, Deligne's work \cite{De10} (and a fortiori Brown's \cite{Br12}) do not predict the nature of $\ast$.

\underline{Case II}: $k>1$ and $n\neq k+1$ (so $n\geq k+3$). Then both quotients in Eq. \eqref{eq30} are singletons. Thus up to isomorphism, there is a unique motive $Z_{n,k}$ with large $\underline{\fu}$ and associated graded isomorphic to $\QQ(n)\oplus \QQ(k) \oplus \mathbbm{1}$. The subobject $W_{-2k} Z_{n,k}$ (resp. subquotient $Z_{n,k}/W_{-2k-1}Z_{n,k}$) of $Z_{n,k}$ is a nontrivial extension of $\QQ(k)$ by $\QQ(n)$ (resp. $\mathbbm{1}$ by $\QQ(k)$). The matrix of periods of $Z_{n,k}$ with respect to suitably chosen bases is of the form
\[
\begin{pmatrix}
(2\pi i)^{-n}&(2\pi i)^{-n}\zeta(n-k) &\ast\\
0&(2\pi i)^{-k}&(2\pi i)^{-k}\zeta(k) \\
0&0&1
\end{pmatrix}.
\]
The period conjecture predicts that $2\pi i, \zeta(k), \zeta(n-k)$ and the entry denoted by $\ast$ are algebraically independent over $\QQ$. Again it would be interesting to find what the entry $\ast$ is. Note that the motive $Z_{n,k}$ is in the subcategory $\mathbf{MT}(\ZZ)$, as from the beginning we may have done the entire discussion of this case in $\mathbf{MT}(\ZZ)$ (as the relevant Ext groups in this case are the same in $\mathbf{MT}(\ZZ)$ and $\mathbf{MT}(\QQ)$). Thus by Brown's work \cite{Br12}, all periods of $Z_{n,k}$ will be in the algebra generated by $2\pi i$ and the multiple zeta values.

\underline{Case III}: $k>1$ and $n=k+1$. This case is the dual situation to Case I. Here the second factor of Eq. \eqref{eq30} is  a singleton, and the motives under investigation are classified up to isomorphism by $Aut(\mathbbm{1})$-orbits of $Ext^1(\mathbbm{1},\QQ(1))$. Consider the complete inequivalent set of representatives $\{[r]\}$ for these orbits as in Case I. Then for each $r$, we get an object $M'_{n,r}$ corresponding to the element of Eq. \eqref{eq30} with the orbit of $[r]$ as its first coordinate. The motives $M'_{n,r}$ are non-isomorphic and up to isomorphism, give all motives with large $\underline{\fu}$ and associated graded isomorphic to $\QQ(n)\oplus \QQ(n-1)\oplus \mathbbm{1}$. 

The motives obtained in this case are intimately related to the $M_{n,r}$ of Case I. Indeed, $M'_{n,r}\dual(n)$ has a large $\underline{\fu}$ (as the property of having a large $\underline{\fu}$ is invariant under dualizing and tensoring by $\QQ(1)$), and its associated graded is isomorphic to $\QQ(n)\oplus \QQ(1)\oplus \mathbbm{1}$. Moreover, the quotient $M'_{n,r}\dual(n)/W_{-2n}$ is isomorphic to the 1-motive $K_r$ given in Eq. \eqref{eq31} (as by construction we have $W_{-2k}M'_{n,r}\simeq K_r(k)$, and $K_r$ is isomorphic to its Cartier dual $K_r\dual(1)$). It follows that $M'_{n,r}\dual(n)$ is isomorphic to $M_{n,r}$ (as they both correspond to the same equivalence class of compatible pairs).

\subsection{}\label{sec: 4 dimensional examples} Let us continue to take $\mathbf{T}=\mathbf{MT}(\QQ)$. The motives of Section \ref{sec: classification of 3 dimensional mixed Tate motives with large u} together with the earlier results of the paper can be used to obtain 4-dimensional mixed Tate motives with 4 weights and a large $\underline{\fu}$.\footnote{Inductively, one can obtain motives with more and more weights which have a large $\underline{\fu}$.} We illustrate this with an example. Let $M$ be the motive $M_{4,r}$ of the previous section, which has associated graded isomorphic to $\QQ(4)\oplus \QQ(1)\oplus \mathbbm{1}$. The weight filtration of $M$ gives an element $\mathcal{L}$ in $Ext^1(\mathbbm{1}, W_{-2}M)$. Let $\mathcal{N}$ a nonzero element of $Ext^1(\mathbbm{1}, \QQ(5))$. Since $Ext^2$ groups vanish in $\mathbf{MT}(\QQ)$, there is an object in $\mathbf{MT}(\QQ)$ attached to the pair 
\[
(\mathcal{L}(5),\mathcal{N}) \ \in \ Ext^1(\QQ(5), (W_{-2}M)(5))\, \times \, Ext^1(\mathbbm{1}, \QQ(5)).
\]
Note that here, at least a priori, there might be non-isomorphic objects attached to the pair, as $Ext^1(\mathbbm{1}, (W_{-2}M)(5))$ is not zero. Any object $\widetilde{M}$ attached to the pair is 4-dimensional, with associated graded isomorphic to
\[
\QQ(9)\oplus \QQ(6)\oplus \QQ(5)\oplus \mathbbm{1}.
\]
Such $\widetilde{M}$ satisfies $(IA3)$, and hence by Theorem \ref{thm3}, $\underline{\fu}(\widetilde{M})$ is large (note that both $M$ and $\mathcal{N}$ have a large $\underline{\fu}$). The field generated over $\QQ$ by the periods of $\widetilde{M}$ contains $2\pi i, \zeta(3), \log r$, the ``new period" of $M$, and $\zeta(5)$. In fact, by the classification of Section \ref{sec: classification of 3 dimensional mixed Tate motives with large u}, the quotient $\widetilde{M}/\QQ(9)$ (which is easily seen to also have a large $\underline{\fu}$) must be isomorphic to the motive $M'_{6,r}$ (of Case III of Section \ref{sec: classification of 3 dimensional mixed Tate motives with large u}), so that the new period of $M'_{6,r}$ will also be a period of $\widetilde{M}$. The period conjecture predicts that the field generated over $\QQ$ by the periods of $\widetilde{M}$ should be of transcendence degree 7 (= ${4\choose 2}+1$), so that $\widetilde{M}$ should have one more new period, which together with the aforementioned six numbers should form an algebraically independent set over $\QQ$.

\begin{rem} 
Note that $k=5$ is the smallest positive integer such that 
\[Gr^WM(k) \oplus \mathbbm{1}\] 
satisfies the independence axiom required to be able to use Theorem \ref{thm3}. 
\end{rem}

\subsection{}\label{counterexample to Thm3 without IA hypothesis} Hypothesis (iii) of Theorem \ref{thm3} was used in the proof to conclude that $\mathcal{E}_p(M)/\underline{\fu}_p(M)$ splits. This hypothesis is actually important for the statement of the theorem to remain true. A counter-example to the statement without this condition can be given in the category $\mathbf{MHS}$ of rational mixed Hodge structures using the work \cite{JR86} of Jacquinot and Ribet on deficient (in the sense of {\it loc. cit.}) points on semiabelian varieties, as we shall discuss below. We shall freely use the basics of the theory of 1-motives (including the realizations of a 1-motive), as introduced by Deligne in \cite{De74}. 

Consider a tuple $(K, A, v, f)$, where
\begin{itemize}
\item[-] $K$ is a number field,
\item[-] $A$ is a simple abelian variety over $K$ with $rank(A(K))>0$,
\item[-] $v\in A^t(K)$ (where $A^t$ is the dual abelian variety),
\item[-] and $f:A^t\longrightarrow A$ is an isogeny over $K$,
\end{itemize}
such that  $f(v)-f^t(v) \, \in \, A(K)$ is a point of infinite order.\footnote{For instance, take $A=A^t$ to be an elliptic curve with complex multiplication by $\ZZ[i]$, $K$ large enough so that complex multiplication by $i$ is defined over $K$ and $A(K)$ has positive rank, $v$ a point of infinite order in $A^t(K)$, and $f=i$ (so that $f^t=-i$).} Let $V$ be a semiabelian variety over $K$, an extension of $A$ by $\mathbb{G}_m$, which under the canonical isomorphism
\[
Ext(A,\mathbb{G}_m) \ \cong \ A^t
\]
corresponds to $v\in A^t(K)$. Denote the projection map $V\longrightarrow A$ by $\pi$. In \cite[Section 4]{JR86}, a point $x_f\in V(K)$ is constructed such that 
\begin{itemize}
\item[(i)] $\pi(x_f)=f(v)-f^t(v)$, and
\item[(ii)] for every nonzero integer $n$ the point $x_f$ is divisible by $n$ in $V(K_n)$, where $K_n$ is the field obtained from $K$ by adjoining the $n$-torsion subgroup of $V$ (such a point is called a deficient point in \cite{JR86}). 
\end{itemize}
Let $M$ be the 1-motive $[\ZZ \, \stackrel{1\mapsto x_f}{\longrightarrow} \, V]$ over $K$. Fixing an embedding $K\subset \overline{K}\subset \CC$, denote the Hodge realization of any 1-motive $N$ over $K$ by $TN$. Thus $TM$ has weights -2,-1, 0 and 
\[
W_{-2}TM \ = \ H_1(\mathbb{G}_m) \ \simeq \ \QQ(1) \, , \quad W_{-1}TM \ = \ H_1(V) \, , \quad Gr^W_{0}TM \ = \ \mathbbm{1}.
\]
We shall see that (with $\mathbf{T}=\mathbf{MHS}$) $\underline{\fu}(TM)$ is not large, whereas both $\underline{\fu}(W_{-1}TM)$ and $\underline{\fu}(TM/W_{-2}TM)$ are large. This would provide a counter-example to the statement of Theorem \ref{thm3} with Hypothesis (iii) of the theorem omitted.

First, let us consider $W_{-1}TM$ and $TM/W_{-2}TM$. The former is a nonsplit extension of the simple Hodge structure $H_1(A)$ by $\QQ(1)$ (because $v$ has infinite order), and hence (by a similar argument as in Section \ref{par: case of two weights}) has a large $\underline{\fu}$. The latter is the Hodge realization of the 1-motive 
\[
[\ZZ \, \stackrel{1\mapsto \pi(x_f)}{\longrightarrow} \, A].
\]
Since $\pi(x_f)$ is a point of infinite order, $TM/W_{-2}TM$ is a nonsplit extension of $\mathbbm{1}$ by $H_1(A)$, and hence has a large $\underline{\fu}$.

To see that $\underline{\fu}(TM)$ is not large, let $\ell$ be a prime number. Given any 1-motive $N$ over $K$, denote the $\ell$-adic realization of $N$ by $T_\ell N$, and let $\Pi_\ell(N)$ be the image of the natural map $Gal(\overline{K}/K)\longrightarrow GL(T_\ell N)(\QQ_\ell)$. Then Property (ii) above implies that the natural (restriction) map
\begin{equation}\label{eq32}
\Pi_\ell(M) \ \longrightarrow \ \Pi_\ell(W_{-1}M)
\end{equation}
(where $W_{-1}M=[0\longrightarrow V])$) is injective (as well as surjective). By the Mumford-Tate conjecture for 1-motives on the unipotent parts (proved by Jossen \cite[Theorem 1]{Jo14}), the Hodge theoretic analogue of this map, i.e. the restriction map
\[
\mathcal{G}(TM,\omega_B) \ \longrightarrow \ \mathcal{G}(T(W_{-1}M),\omega_B)\hspace{.4in}(\text{$\omega_B$ = the forgetful fiber functor})
\]
is also injective (the two groups above are calculated in $\mathbf{MHS}$). Thus $\underline{\fu}_{-1}(TM)$ is zero.

\begin{rem} \ \\
\vspace{-.1in}
\begin{itemize}
\item[(1)] Here we do not need the full power of the Mumford-Tate conjecture on the unipotent parts to go from injectivity of Eq. \eqref{eq32} to the vanishing of $\underline{\fu}_{-1}(TM)$; just the more basic statement \cite[Theorem 1]{Ber98} is enough. Indeed, \cite[Theorem 1]{Ber98} and injectivity of Eq. \eqref{eq32} imply that $W_{-2}\underline{\fu}(TM)$ is zero. It follows that $\underline{\fu}(TM)$ and consequently $\underline{\fu}_{-1}(TM)$ is a pure object of weight -1. On the other hand, 
\[
\underline{\fu}_{-1}(TM)\subset \inHom(TM/W_{-1}TM, W_{-1}TM)\cong W_{-1}TM.
\]
It follows that $\underline{\fu}_{-1}(TM)$ is zero (as otherwise, in light of simplicity of $H_1(A)$ the extension $\mathcal{E}_{-2}(W_{-1}TM)$ would split).

\item[(2)] Note that the example given in this section shows that in general, without any independence axiom, the individual extensions $\mathcal{E}_p/\underline{\fu}$ need not split (see Corollaries \ref{cor 1 of thm 2} and \ref{cor 2 of thm 2} of Theorem \ref{thm 2}, as well as Deligne's Theorem \ref{Deligne thm} and the Remark after). Indeed, in the above example, $\mathcal{E}_{-1}(TM)/\underline{\fu}(TM)$ does not split: if it did, then by Lemma \ref{lem equivalent statements to E_p/u_p originaring from S} so would $\mathcal{E}_{-1}(TM)/\underline{\fu}_{-1}(TM)$. But $\mathcal{E}_{-1}(TM)$ ( = $\mathcal{E}_{-1}(TM)/\underline{\fu}_{-1}(TM)$) does not split as $x_f$ is a point of infinite order.

\item[(3)] In fact, the example given in this section also shows that in general, $\underline{\fu}$ may not be the sum of the subobjects $\underline{\fu}_p$. Indeed, with $M$ as above, $\underline{\fu}(TM)$ is not zero (because $TM$ is not semisimple), while both $\underline{\fu}_{-1}(TM)$ and $\underline{\fu}_{-2}(TM)$ are zero. (That the latter is zero can be seen by an argument similar to the one given in Remark (1): $\underline{\fu}_{-2}(TM)$ is pure of weight -1 and a subobject of $\inHom(TM/W_{-2}TM,W_{-2}TM)\cong  (TM/W_{-2}TM)\dual(1)$; the latter object has no nonzero subobject of weight -1.)

\end{itemize}
\end{rem}


\begin{thebibliography}{00}


%
\bibitem{An04}
Y. Andr\'{e},
Une introduction aux motifs,
Soci\'{e}t\'{e} Math\'{e}matique de France, 2004

\bibitem{Ar13}
D. Arapura, 
An abelian category of motivic sheaves, 
Advances in Math. 233 1 (2013), 135-195.

\bibitem{Ay14}
J. Ayoub, 
L’alg\`{e}bre de Hopf et le groupe de Galois motiviques d’un corps de caract\'{e}ristique nulle, 
J. Reine Angew. Math. 693 (2014), partie I: 1-149 ; partie II: 151-226

\bibitem{Be83}
A. A. Beilinson, 
Notes on absolute Hodge cohomology, 
Applications of Algebraic K-theory to Algebraic Geometry and Number Theory, Part I, Proceedings of a Summer Research Conference held June 12-18, 1983, in Boulder, Colorado, Contemporary Mathematics 55, American Mathematical Society, Providence, Rhode Island, pp. 35–68

\bibitem{Be02}
C. Bertolin,
The Mumford-Tate group of 1-motives,
Ann. Inst. Fourier, Grenoble, 52, 4 (2002), 1041-1059

\bibitem{Be03}
C. Bertolin,
Le radical unipotent du groupe de Galois motovique d'un 1-motif,
Math. Ann. 327, 585-607 (2003)

\bibitem{Ber01}
D. Bertrand,
Unipotent radicals of differential Galois group and integrals of solutions of inhomogeneous equations,
Math. Ann. 321 (2001), no. 3, 645–666

\bibitem{Ber98}
D. Bertrand, 
Relative splitting of 1-motives,
Contemporary Mathematics, Vol. 210, 1998

\bibitem{Ber13}
D. Bertrand,
Extensions panach\'{e}es autoduales,
J. K-Theory 11 (2013), no. 2, 393–411


\bibitem{Br12}
F. Brown,
Mixed Tate motives over $\ZZ$,
Annals of Mathematic (175), 2012, 949-976


\bibitem{Br16}
F. Brown,
Irrationality proofs for zeta values, moduli spaces and dinner parties,
Mosc. J. Comb. Number Theory 6 (2016), 102-165


\bibitem{Bu59}
D. A. Buchsbaum,
A note on homology in categories,
Annals of Mathematics, Second Series, Vol. 69, No. 1 (Jan., 1959), pp. 66-74 


\bibitem{Ca80}
J. A. Carlson,
Extensions of mixed Hodge structures,
Journ\'{e}es de G\'{e}ometrie Alg\'{e}brique d'Angers, Juillet 1979/Algebraic Geometry, Angers, 1979, pp. 107–127, Sijthoff \& Noordhoff, Alphen aan den Rijn, Md., 1980.



\bibitem{De74}
P. Deligne,
Theorie de Hodge III,
Publications Math\'{e}matiques de l'I.H.\'{E}.S., tome 44 (1974), p. 5-77



\bibitem{De10}
P. Deligne,
Le groupe fondamental motivique de $G_m-\mu_N$ pour $N=2,3,4,6$ ou $8$.
Publications Math\'{e}matiques de l'IH\'{E}S. 112 (2010) pp. 101-141

\bibitem{DG05}
P. Deligne and A. B. Goncharov,
Groupes fondamentaux motiviques de Tate mixte,
{\it Ann. Sci. de l’\'{E}cole Norm. Sup.}
38 (2005) 1-56 

\bibitem{DM82}
P. Deligne and J.S. Milne,
Tannakian Categories,
In Hodge Cycles, Motives, and Shimura Varieties,
Lecture Notes in Mathematics 900, Springer-Verlog, Berlin (1982)

\bibitem{De89}
P. Deligne,
Le groupe fondamental de la droite projective moins trois points,
Galois Groups over $\QQ$, MSRI Publ. 16, pp. 79-313, Springer-Verlag (1989)

\bibitem{De90}
P. Deligne,
Cat\'{e}gories tannakiennes,
Grothendieck Festschrift Vol II, Progress in Mathematics, 87 ( Birkhäuser Boston 1990) pp. 111–195

\bibitem{De94}
P. Deligne,
Structures de Hodge mixtes r\'{e}elles, 
Proceedings of Symposia in Pure Mathematics, Volume 55 (1994), Part I



\bibitem{DW16}
I. Dan-Cohen and S. Wewers, 
Mixed Tate motives and the unit equation,
Int. Math. Res. Not., 2016, no. 17.

\bibitem{Dunham}
W. Dunham, Euler and the Cubic Basel Problem, {\em The American Mathematical Monthly}, 128(2021), 291-301.


\bibitem{EM21}
P. Eskandari and V. K. Murty,
The fundamental group of an extension in a Tannakian category and the unipotent radical of the Mumford-Tate group of an open curve, preprint, available on \url{http://www.math.utoronto.ca/payman}

%

\bibitem{Euler}
L. Euler, De relatione inter ternas pluresve quantitates instituenda (E591), Opusc. Anal. 2(1785), 91?101. 
Reprinted in Opera Omnia. Ser. 1, Vol. 4: 136-145.

\bibitem{Go05}
A. B. Goncharov,
Galois symmetries of fundamental groupoids and noncommutative geometry,
Duke Math. Journal, 128, no. 2 (2005) p. 209-284


\bibitem{Gr68}
A. Grothendieck,
Mod\`{e}les de N\'{e}ron et monodromie,
SGA VII.1, no 9, Springer LN 288, 1968




%

\bibitem{Har06}
C. Hardouin,
Hypertranscendance et Groupes de Galois aux diff\'{e}rences, arXiv 0609646v2, 2006

\bibitem{Har11}
C. Hardouin,
Unipotent radicals of Tannakian Galois groups in positive characteristic, 
Arithmetic and Galois theories of differential equations, 223-239, S\'{e}min. Congr., 23, Soc. Math. France, Paris, 2011

\bibitem{HM17}
A. Huber and S. M\"{u}ller-Stach,
Periods and Nori Motives,
Springer, 2017

\bibitem{JR86}
O. Jacquinot and K. Ribet,
Deficient points on extensions of abelian varieties by G,
J. Number Theory, 24, No. 3 (1986).

\bibitem{Ja90}
U. Jannsen,
Mixed motives and algebraic K-theory,
with appendices by S. Bloch and C. Schoen,
Lecture Notes in Mathematics, Vol. 1400,
Springer-Verlag, Berlin, 1990.

\bibitem{Ja94}
U. Jannsen,
Motivic sheaves and filtrations on Chow groups,
Proceedings of Symposia in Pure Mathematics,
Volume 55 (1994), Part 1

\bibitem{Jo14}
P. Jossen,
On the Mumford-Tate conjecture for 1-motives,
Inventiones Math. (2014) 195: 393-439



\bibitem{Mi17}
J. S. Milne,
Algebraic Groups,
Cambridge University Press, 2017

\bibitem{Ne94}
J. Nekovar, 
Beilinson's conjectures, 
in Motives (Seattle, WA, 1991), 537 - 570, Proc. Symp. Pure Math., 55/I, Amer. Math. Soc., Providence, RI, 1994




\bibitem{Riv72}
N. Saavedra Rivano.
Cat\'{e}gories Tannakiennes,
Lecture Notes in Mathematics 265

%

\bibitem{Yo60}
N. Yoneda,
On Ext and exact sequences,
Journal of the Faculty of Science, Imperial University of Tokyo, Vol. 8, 1960,  507-576

\end{thebibliography}
\end{document}